\documentclass[12pt,leqno,titlepage]{article}

\usepackage{amsmath}
\usepackage{amsfonts}
\usepackage{amssymb}
\usepackage{euscript}

\usepackage{graphicx}
\usepackage{latexsym} 

\usepackage{accents}

\usepackage{xcolor}

\usepackage{hyperref}
\hypersetup{
 bookmarks=true,
 pdftitle={Meanders and hemispheres},
 pdfauthor={Carlos Rocha, Bernold Fiedler},
 colorlinks=true,
 linkcolor=blue,
 citecolor=blue,
 filecolor=blue,
 urlcolor=blue}

\newtheorem{theorem}{Theorem}[section] 
\newtheorem{lemma}{Lemma}[section] \baselineskip 1.5em

\numberwithin{equation}{section}
\numberwithin{figure}{section}

\providecommand{\N}{\mathbb{N}}
\providecommand{\Z}{\mathbb{Z}}
\providecommand{\R}{\mathbb{R}}

\newcommand{\sign}{\operatorname{sign}}

\def\cE{\EuScript{E}}
\def\cA{\EuScript{A}}
\def\cC{\EuScript{C}}
\def\cM{\EuScript{M}}
\def\cO{\EuScript{O}}
\def\cJ{\EuScript{J}}
\def\cH{\EuScript{H}}

\title{{\Large Meanders, zero numbers and the \\ cell structure of Sturm global attractors}\\
\vspace{1cm}
{\vspace{1ex}{\large -- Dedicated to the memory of Pavol Brunovsk\'y --}}\vspace{1ex}}

\author{Carlos Rocha\\
        Instituto Superior T\'ecnico\\
        Avenida Rovisco Pais, 1049--001 Lisboa, PORTUGAL\\
        {\tt crocha@math.ist.utl.pt}\\
        {\tt http://www.math.ist.utl.pt/cam/}\\
        {}
  \and
        Bernold Fiedler\\
        Institut f\"ur Mathematik\\
        Freie Universit\"at Berlin\\
        Arnimallee 7, D--14195 Berlin, GERMANY\\
        {\tt fiedler@math.fu-berlin.de}\\
        {\tt http://dynamics.mi.fu-berlin.de}\\
        {}  }

\date{version of \today}

\begin{document}

\maketitle

\setlength{\parindent}{0cm}

\begin{abstract}
We study global attractors $\cA=\cA_f$ of semiflows generated by semilinear partial 
parabolic differential equations of the form $u_t = u_{xx} + f(x,u,u_x), 0<x<1$, satisfying 
Neumann boundary conditions. The equilibria $v\in\cE\subset\cA$ of the semiflow are 
the stationary solutions of the PDE, hence they are solutions of the corresponding 
second order ODE boundary value problem. 
Assuming hyperbolicity of all equilibria, the dynamic decomposition 
of $\cA$ into unstable manifolds of equilibria provides a geometric and topological 
characterization of Sturm global attractors $\cA$ as finite regular signed CW-complexes,
the Sturm complexes, with cells given by the unstable manifolds of equilibria.
Concurrently, the permutation $\sigma=\sigma_f$ derived from the ODE boundary value 
problem by ordering the equilibria according to their values at the boundaries $x=0,1$, 
respectively, completely determines the Sturm global attractor $\cA$. 
Equivalently, we use a planar curve, the meander $\cM=\cM_f$, associated to the the 
ODE boundary value problem by shooting. 

In the previous paper \cite{firo18d}, we set up to determine the boundary neighbors of any specific 
unstable equilibrium $\cO$, based exclusively on the information on the corresponding signed 
hemisphere complex. In addition, a certain minimax property of the boundary neighbors was established. 
In the signed hemisphere decomposition of the cell boundary of $\cO$,
this property identifies the equilibria which are closest to, or most 
distant from, $\cO$ at the boundaries $x=0,1$, in each hemisphere. 

The main objective of the present paper is to derive this minimax property directly from the Sturm 
permutation $\sigma$, or equivalently from the Sturm meander $\cM$, based 
on the Sturm nodal properties of the solutions of the ODE boundary value problem. 
This minimax result simplifies the task of identifying the equilibria on the cell boundary of each 
unstable equilibrium, directly from the Sturm meander $\cM$.
We emphasize the local aspect of this result by an example for which the identification 
of the equilibria is obtained from the knowledge of only a segment of 
the Sturm meander $\cM$.

\end{abstract}

\section{Introduction}\label{sec1}

For quite a long time we have been pursuing the study of global attractors of scalar semilinear 
parabolic equations in one space dimension. 
Under appropriate dissipativeness conditions on the nonlinearity $f\in C^2$, equations of the form
\begin{equation}\label{eq101} 
u_t = u_{xx} + f(x,u,u_x) \,, \ 0 < x < 1 \,,
\end{equation}
with suitable boundary conditions, generate global dissipative semiflows on a Sobolev space 
$X \subset C^1([0,1],\R)$. The global attractor $\cA_f\subset X$ of the semiflow, i.e.~the maximal 
compact invariant subset, is called {\em Sturm global attractor}. 
In the present paper we consider the specific case of Neumann boundary conditions $u_x=0$ at $x=0,1$, 
for which the generated semiflow is gradient-like. 
We refer to \cite{hen81,paz83,tan79} for a general background, to \cite{hal88,tem88,lad91,bavi92,rau02} 
for a contextual introduction to global attractors, and to \cite{firo14,firo15}, and the citations 
there, for a general introduction to Sturm global attractors. See also \cite{brfi88,brfi89,fie94} 
for earlier results strongly influenced by Palo Brunovsk\'y.

The stationary solutions $u(t,x)=v(x)$ of \eqref{eq101} satisfy the second order Neumann boundary 
value problem
\begin{equation}\label{eq102} 
\begin{aligned}
v_{xx} + f(x,v,v_x) & = 0 \,, \ 0 < x < 1 \,, \\
v_x & = 0 \,, \ x=0 \mbox{ or } 1 \,,
\end{aligned}
\end{equation}
which provides the link to the name of Sturm associated to the global attractor $\cA=\cA_f$. 
Each stationary solution of \eqref{eq101} corresponds to an equilibrium $v\in\cE$ of the 
semiflow in $X$ and is a compact invariant subset, therefore an element of the global attractor $\cA$. 
Under the generic assumption of {\em hyperbolicity} of all equilibria $v\in\cE$ the set of 
equilibria is finite, 
\begin{equation}\label{eq103} 
\cE = \{v_j : j=1,\dots,N \} \,,
\end{equation}
with $N$ odd due to dissipativeness. Hyperbolicity of the equilibrium $v\in\cE$ is achieved if  
$\lambda=0$ is not an eigenvalue of the linearization of \eqref{eq101} at  $v$.
Then, the gradient-like property of the semiflow provides a dynamic decomposition of the Sturm global 
attractor $\cA$ into unstable manifolds of equilibria,
\begin{equation}\label{eq104} 
\cA = \bigcup_{1\le j\le N} W^u(v_j) \,.
\end{equation}
This decomposition suggests a geometric and topological characterization of the Sturm global 
attractor $\cA$ as a finite $CW$-complex $\cC=\cC_f$,  
\begin{equation}\label{eq105} 
\cC = \bigcup_{1\le j\le N} c_j \,,
\end{equation}
with cell interiors $c_j=W^u(v_j)$. 
Although this idea fails, in general gradient settings, it works well in our Sturm setting \eqref{eq101}.
In fact we obtain the regular {\em Thom-Smale complex} of the Sturm global 
attractor $\cA$, or the {\em Sturm complex}, \cite{firo14,firo18a}.

Concurrently, the permutation $\sigma=\sigma_f\in S_N$ derived from the boundary value problem 
\eqref{eq102} by ordering successively the $N$ equilibria according to their values at $x=0$ and $x=1$, 
completely determines the Sturm global attractor $\cA$ up to $C^0$ orbit equivalence, \cite{firo00}.
If we let $h_0,h_1: \{1,\dots,N\}\rightarrow\cE$ denote the two boundary orders of the equilibria,
\begin{equation}\label{eq106} 
h_\iota(1)<h_\iota(2)<\dots<h_\iota(N) \quad \mbox{ at } \quad x=\iota\in\{0,1\} \,,
\end{equation}
then the {\em Sturm permutation} 
\begin{equation}\label{eq107} 
\sigma := h_0^{-1}\circ h_1
\end{equation}
provides a combinatorial description of the Sturm global attractor $\cA$ determining exactly 
which equilibria $v_j,v_k\in\cE$ are connected by heteroclinic connecting orbits, $v_j\leadsto v_k$. 
For this combinatorial description of Sturm global attractors we refer to \cite{furo91,firo96}.
Here we just mention the central role of the nodal properties of the solutions of \eqref{eq101}. 
If $0\le z(\varphi)\le\infty$ denotes the number of strict sign changes (the {\em zero number}) of 
$\varphi:[0,1]\rightarrow\R$, with $\varphi\not\equiv 0$, then 
\begin{equation}\label{eq108}
t \quad \longmapsto \quad z(u^1(t,\cdot)-u^2(t,\cdot))
\end{equation}
is finite and nonincreasing for $t>0$, for any two distinct solutions $u^1,u^2$ of \eqref{eq101},
and drops strictly whenever multiple zeros $u^1 = u^2, u^1_x = u^2_x$ occur at any $t_0, x_0$. 
We refer to \cite{ang88} for details and to \cite{mat82,firo96,firo99,firo00,roc91,gal04} 
for many aspects of nonlinear Sturm theory. 
For convenience, in the following we also use the signed notation 
\begin{equation}\label{eq109}
z(\varphi) = z_\pm
\end{equation}
to indicate that $\varphi$ has $z$ strict sign changes and the index $\pm$ to indicate 
that $\pm\varphi(0)>0$.

The solution of the second order boundary value problem \eqref{eq102} by ODE shooting in phase space 
$(v,v_x)$ produces a meander characterization of the Sturm permutation. 
Existence of solutions for $x\in[0,1]$ is ensured here by assuming a sublinear growth of $f$, without 
loss of generality for $\cA$ due to its compactness. 
Let $v=v(x,a)$ abbreviate the solution starting from $v=a, v_x=0$ at $x=0$.
Then, if we consider the initial set of Neumann conditions $\{(a,0): a\in\R\}$ (at $x=0$) we obtain as 
image at $x=1$ a planar regular non-selfintersecting curve $\cM=\cM_f:=\{(v(1,a),v_x(1,a)): a\in\R\}$ 
which transversely intersects the Neumann horizontal axis at $N$ points corresponding to the end values 
of the equilibria $v_j\in\cE$. 
Such a curve is called a {\em meander} and its existence shows that the Sturm permutation is 
a {\em meander permutation}. 
We refer to \cite{firo99} for details, and to \cite{kar17,detal19} for many other aspects of meander 
permutations.
The first important outcome of the meander characterization is that the ODE Sturm permutation 
$\sigma=\sigma_f$ determines the PDE {\em Morse indices} $i(v_j)=\dim W^u(v_j)$ of the equilibria 
$v_j\in\cE, j=1,\dots,N$, and the number of zeros of their differences $z(v_k-v_j), 1\le j<k\le N$.  
We review these results in sections \ref{sec2} and \ref{sec5} below. For now it suffices to recall that, 
for odd $N$, a permutation $\sigma\in S_N$ is a Sturm permutation if and only if it is a 
{\em dissipative}, {\em Morse}, meander permutation, \cite{firo99}. 
Dissipative here means that $\sigma(1)=1$ and $\sigma(N)=N$ are fixed, and Morse means that the Morse 
indices computed from $\sigma$ are all non-negative, $i(v_j)\ge 0$.


Let $\cO\in\cE$ denote an unstable equilibrium with $i(\cO)=n$, and let $\cE^k_\pm(\cO)\subset\cE$ 
denote the equilibrium subsets
\begin{equation}\label{eq113} 
\cE^k_\pm(\cO) := \{w\in\cE: z(w-\cO)=k_\pm \mbox{ and } \cO\leadsto w \} \,.
\end{equation} 
The four {\em boundary neighbors} $w^\iota_\pm = w^\iota_\pm(\cO)$ of $\cO$,
\begin{equation} \label{eq111} 
w^\iota_\pm := h_\iota(h_\iota^{-1}(\cO)\pm1) \quad  , \quad \iota\in\{0,1\} \,,
\end{equation}
whenever defined, are the predecessors and successors of $\cO$, along the boundary orders $h_\iota$ at 
$x=\iota\in\{0,1\}$.
By Lemma \ref{l52} below, any boundary neighbor $w$ of $\mathcal{O}$, separately, either possesses Morse 
index $i(w)=n-1$ or else Morse index $i(w)=n+1$.
In \cite{firo18d}, the identification of the boundary neighbors of $\cO\in\cE$ with 
\begin{equation} \label{eq111a} 
i(w^\iota_\pm) = i(\cO)-1 = n-1 \,,
\end{equation} 
was obtained from the signed zero numbers \eqref{eq109} of the differences $\varphi=w-\cO$ for all 
$w\in\cE^k_\pm(\cO)$, using the zero number decay property \eqref{eq108}. 
Quite surprisingly it turned out that the 
closest equilibrium to $\cO$ at the boundary $x=\iota$ is also the most distant to $\cO$ in the 
opposite boundary $x=1-\iota$.
We call this feature the {\it minimax property}. 

Restricting our attention to the equilibria in $\cE^{n-1}_\pm(\cO)$, we let 
$\underline v^{\,\iota}_{\,\pm}, \overline v_{\,\pm}^{\,\iota}$, $\iota\in\{0,1\}$, denote the 
{\it minimax equilibria} defined by  
\begin{eqnarray} \label{eq114}
\underline v^{\,\iota}_{\,\pm} & := & \mbox{ the equilibrium } v\in\cE_\pm^{n-1}(\cO) \mbox{ closest to }
\cO \mbox{ at } x=\iota \,, \\ \label{eq115}
\overline v_{\,\pm}^{\,\iota} & := & \mbox{ the equilibrium } v\in\cE_\pm^{n-1}(\cO) \mbox{ most distant 
from } \cO \mbox{ at } x=\iota \,. 
\end{eqnarray}

For example, let us restrict attention to those boundary neighbors $w=w^\iota_\pm$ of $\cO$ defined 
in \eqref{eq111} which satisfy $i(w)=n-1$, if any.
In section \ref{sec2} we show that $w\in \cE^{n-1}_\pm(\cO)$. 
In fact we prove, much more precisely:
%
%
\begin{eqnarray}
	w^0_-&=& \underline v_{\,-}^{\,0}; \label{eq116}\\
	w^0_+&=& \underline v^{\,0}_{\,+}; \label{eq116a}\\
	w^1_-&=& 
	   \begin{cases}
	   \underline v^{\,1}_{\,+}, \quad \text{for even } n;\\
     \underline v_{\,-}^{\,1}, \quad \text{for  odd } n;
	   \end{cases} \label{eq116b}\\
	w^1_+&=& 
	   \begin{cases}
	   \underline v_{\,-}^{\,1}, \quad \text{for even } n;\\
     \underline v^{\,1}_{\,+}, \quad \text{for  odd } n.\\
	   \end{cases} \label{eq116c}
\end{eqnarray}

The surprising \emph{minimax property} now asserts the equality of the equilibria defined by \eqref{eq114} 
and \eqref{eq115}: 
\begin{equation} \label{eq117} 
\underline v^{\,\iota}_{\,\pm} = \overline v_{\,\pm}^{\,1-\iota} \,.
\end{equation} 

This provides equivalent expressions for the meander neighbors $w$ of $\cO$ in terms of the 
most distant equilibria $\overline v_{\,\pm}^{\,1-\iota}$, rather than the closest neighbors 
$\underline v^{\,\iota}_{\,\pm} $ in \eqref{eq116}--\eqref{eq116c}.

Evidently, the minimax property also simplifies the task of identifying the equilibria in 
$\cE^{n-1}_\pm(\cO)$, especially when viewed directly in the Sturm meander $\cM$.
For example, it is interesting to check our results against the influential original description of all 
heteroclinic targets $\cO \leadsto v$, for nonlinearities $f=f(u)$. 
See \cite{brfi88,brfi89}, Theorems 1.4 and 1.5 (where $\cO$ was called $v$ and $v$ was called $w$).
Obtained under Dirichlet conditions, those results have to be adapted slightly to reveal a minimax theorem 
for the ordering of equilibria by boundary derivatives $v_x$, in the target set $\Omega_3$\,, rather than 
by boundary values $v$.
See \cite{firowo12} for a more recent discussion of the Neumann case.

The main purpose of our present paper is to prove the minimax property, Theorem~\ref{theo11}, 
based solely on plane meander arguments of ODE type. 
We believe that, by its simplicity, this more elementary approach further elucidates the structure of 
global Sturm attractors, and is therefore of independent interest. 

\begin{theorem} \label{theo11}
Let $\underline v^{\,\iota}_{\,\pm}, \overline v_{\,\pm}^{\,\iota}\in\cE_\pm^{n-1}(\cO)$ denote the 
minimax equilibria of $\cO$ defined in \eqref{eq114}--\eqref{eq115}. If any of the $w_\pm^\iota$
of the unstable equilibrium $\cO$, as defined by \eqref{eq111}, is more stable than $\cO$, i.e. if 
$i(w_\pm^\iota)=n-1$ as in \eqref{eq111a}, then, for the associated $\underline v^{\,\iota}_{\,\pm}$ of 
that sign $\pm$ and that $\iota$, in \eqref{eq114}--\eqref{eq116c}, 
the minimax property \eqref{eq117} holds. 
\end{theorem}

The minimax property \eqref{eq117} was also obtained in Theorem~4.3 of \cite{firo18d} -- 
regardless of the particular Morse indices $i(w_\pm^\iota)$ of the immediate 
$h_\iota$\nobreakdash-neighbors $w_\pm^\iota$ of $\cO$.
Moreover, the minimax property was proved for all minimax equilibria in the equilibrium sets 
$\cE^k_\pm(\cO),\ 0\leq k < n=i(\cO)$, not just for $k=n-1$.

The proof, however, was based on the \emph{signed hemisphere complex} $\Sigma^k_\pm(\cO)$ associated 
to any Sturm global attractor $\cA$. 
The hemispheres arise from the decompositions of the $k$-dimensional sphere boundaries 
$\Sigma^k=\partial W^{k+1}(\cO)$ of the $(k+1)$-dimensional fast unstable manifolds $W^{k+1}(\cO)$ 
by the $(k-1)$-sphere boundary of the $k$-dimensional faster unstable manifold $W^{k}(\cO)$.
A crucial bridge to the current meander-based definition \eqref{eq113} of $\cE^k_\pm(\cO)$ was then 
the observation
\begin{equation}
\label{eq118}
\cE^k_\pm(\cO) = \cE \cap \Sigma^k_\pm(\cO)\,.
\end{equation}
Finally, the minimax equilibria were identified recursively, from just the incidence geometry of the 
hemisphere cells of adjacent dimensions.
This allowed us to determine the minimax equilibria, and thereby the boundary neighbors $w$ of $\cO$ in case 
$i(w)=n-1$, from purely geometric properties of the given signed hemisphere complex $\Sigma^k_\pm(\cO)$.

Our results therefore gain broader relevance as a step within a much more ambitious program: the complete 
geometric characterization of all those signed regular cell complexes which arise as signed Thom-Smale 
complexes in the Sturm PDE setting. Limited to mere 3-ball global attractors $\mathcal{A}$, we initiated 
this program in \cite{firo18a}--\cite{firo18c}.
Even with the step presented here and in the companion paper \cite{firo18d}, we are still far 
from that goal, in general.

The missing steps would involve the following.
Of course, we first have to show that the (still elusive) characterizing geometric properties of the 
initial formal signed regular cell complex, are satisfied by all Sturm complexes.
Given an abstract cell complex, prescribed according to the rules of the geometric characterization, we 
formally call the barycenters of the abstract cells ``equilibria''.
The cell dimensions are the ``Morse indices'' of the barycenter equilibria.
The rules of Theorem~\ref{theo11} then formally determine predecessors and successors of all ``equilibria'', 
recursively.

Such characterizing rules have only been obtained, successfully, in low dimensional cases of 
\begin{equation}\label{eq110} 
\dim\cA:=\max\{\dim W^u(v_j), 1\le j\le N\} \,. 
\end{equation}
For the planar case see \cite{firo09a,firo08,firo09b}; 
for three-dimensional balls see \cite{firo18a,firo18b,firo18c}.

The results in \cite{firo18d} determine the boundary neighbors \eqref{eq111} of any unstable equilibrium, 
based on the geometry of a signed hemisphere complex, exclusively.
Starting from the ``equilibria'' of maximal ``Morse index'',~i.e. from the barycenters of cells with 
maximal dimension, this recursively determines two formal ``boundary orders'' $h_\iota$ of all 
equilibria, separately for each $\iota=0,1$.
The two formal ``total orders'' $h_\iota$\,, in turn, determine a formal permutation $\sigma$, as in 
\eqref{eq106}, \eqref{eq107}.
The elusive geometric characterization rules, however, would have to guarantee, a priori, that the two 
resulting objects $h_\iota$ define total orders, i.e.~Hamiltonian paths, among the equilibria $\cE$.
Moreover, the formal permutation $\sigma := h_0^{-1}\circ h_1$ has to turn out as a Sturm permutation, in 
particular a meander.
In a final step, it remains to show that the resulting signed Thom-Smale complex, associated to the Sturm 
permutation $\sigma$, coincides with the prescribed original cell complex.
In particular, the barycenter ``equilibria'' become equilibria, their cells are their unstable manifolds, 
the ``Morse indices'' become actual Morse indices, i.e. unstable dimensions, and the formal total 
``boundary'' orders $h_\iota$ become the orders of the equilibria at the respective boundaries $x=\iota=0,1$.
We repeat that this program has only been completed for planar and 3-ball Sturm attractors, so far.
For further illustration in the present context we refer to the discussion section of \cite{firo18d}.

In the above program of a geometric characterization of all Sturmian signed hemisphere complexes, our 
present paper pursues a complementary approach.
Here we characterize boundary neighbors by our minimax theorem \ref{theo11}, based on meander arguments only.
We expect the insights gained here to play an essential role in our attempt to bridge the gap between the 
meander description and the geometric cell complex description of Sturm global attractors, eventually.
Another advantage of our present approach is the relatively elementary ODE level required.

In the next section \ref{sec2} we recall notation and essential results necessary to deal with the 
Morse indices of equilibria and the zero numbers of their differences, directly from the 
meander $\cM$ or its meander permutation. In section \ref{sec3} we prove the main result, 
Theorem~\ref{theo11}.
To emphasize a local aspect of our global result, in section \ref{sec4} we show an example for which 
the identification of the equilibria in $\cE^{i(\cO)-1}_+(\cO)$ is obtained from the knowledge of 
only a section of the meander $\cM$. 
We finish with an appendix, section \ref{sec5}, where we review the nonlinear Sturm-Liouville property in 
our meander setting, and another appendix, section \ref{sec6}, where we describe the double cone suspension 
of Sturm global attractors used in the proof of Theorem~\ref{theo11}.

\vspace*{.3cm}

\textbf{Acknowledgments.} 
This paper is dedicated to the memory of our longtime friend, colleague, and coauthor Palo Brunovsk\'y 
in deep admiration and gratitude. We owe much of our quest to his pioniering curiosity and his friendly, 
sharing, and noble spirit.
Extended mutually delightful hospitality by the authors has gratefully been enjoyed.
CR expresses also gratitude to his family, friends and longtime coauthor in appreciation of their support 
and patience during a recent and specially hard time.
This work was partially supported by DFG/Germany through SFB 910 project A4, and by FCT/Portugal through 
projects UID/MAT/04459/2019 and UIDB/04459/2020.

\section{Meanders and Crossing numbers} \label{sec2}

\begin{figure}[t] 
\centering \includegraphics[scale = 1.0]{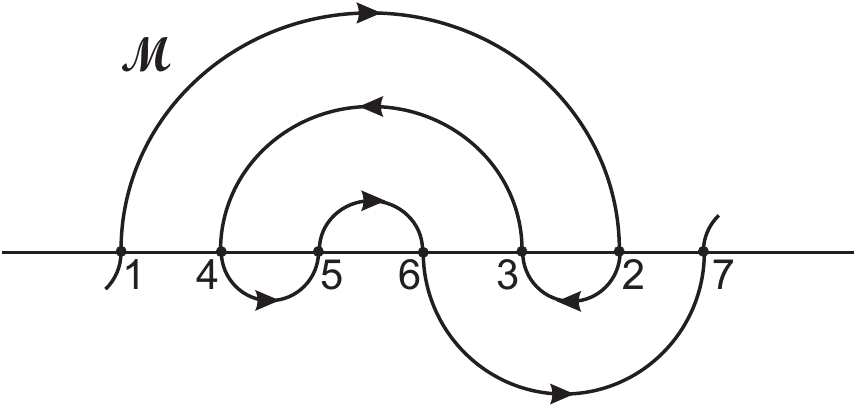}
\caption{\em\small A meander $\cM$ in canonical form. The intersection points are labeled along the meander.
Along the horizontal axis the labeling corresponds to the permutation 
$\sigma = \{ 1 \ 4 \ 5 \ 6 \ 3 \ 2 \ 7 \}$.}
\label{fig21}
\end{figure}

A {\em meander} is a $C^1$-smooth planar embedding of the (oriented) real line which crosses a 
positively oriented horizontal line at finitely many points with transverse intersections, \cite{arn88}.
We consider meanders that run from Southwest to Northeast asymptotically, hence the number $N$ of 
crossing points is odd. A meander for which all intersections with the base line are vertical, 
and all arcs joining intersection points are semicircles, is said to be in {\em canonical form} 
(see Figure~\ref{fig21} for an example).

In the following, we label the intersection points along the meander. The labeling obtained by 
reading along the horizontal axis corresponds to a permutation $\sigma\in S_N$. Any permutation 
$\sigma\in S_N$ associated to a meander is then called a {\em meander permutation}. 
We recall that a permutation $\sigma\in S_N$ is called {\em dissipative} if it fixes the end points, 
i.e.~$\sigma(1)=1$ and $\sigma(N)=N$.
The final ingredient has to do with the winding of the unit tangent vector of the meander. 
Solely based on the permutation $\sigma\in S_N$ we inductively define the Morse numbers $i_j$ of the 
intersection points labeled $j=1,\dots,N$ by
\begin{equation}\label{eq201} 
i_1 = 0 \,, \quad i_{j+1} = i_j + (-1)^{j+1} \sign\left(\sigma^{-1}(j+1)-\sigma^{-1}(j)\right) \,.
\end{equation}
With the meander in canonical form, these numbers count the clockwise half-windings of the unit 
tangent vector of the meander, from the initial intersection point $j=1$ to the intersection point 
$j\in\{1,\dots,N\}$. 
We say that $\sigma\in S_N$ is a {\em Morse permutation} if all its Morse numbers are non-negative, 
that is, $i_j\ge 0$ for all $j=1,\dots,N$.

A dissipative, Morse, meander permutation $\sigma\in S_N$ is called a {\em Sturm permutation}, 
and a meander with an associated Sturm permutation is called a {\em Sturm meander}.
As it follows from \cite[Theorem 1.2]{firo99}, any Sturm permutation is realizable by a boundary 
value problem \eqref{eq102} with dissipative nonlinearity $f$, and the corresponding Sturm meander 
is equivalent to the associated ODE shooting meander in phase space $(u,u_x)\in\R^2$.
Moreover, since global Sturm attractors $\cA_{f_1}$ and $\cA_{f_2}$ with the same Sturm permutation 
$\sigma_{f_1}=\sigma_{f_2}$ are $C^0$ orbit equivalent (see \cite{firo00}), from here on we assume 
all meanders are in canonical form.

Let $\cM$ denote a Sturm meander intersecting the horizontal axis.
Let the intersection points be labeled by $j=1,\dots,N$, enumerated along the meander. 
For any Sturm realization \eqref{eq102} of $\cM$, we label the set of equilibria $\cE=\{v_j, j=1,\dots,N\}$ 
such that $v_j=h_0(j)$.
Then their Morse indices are given by the Morse numbers, $i(v_j)=i_j$. 
The interpretation of Figure~\ref{fig21} as a shooting diagram and the ordering \eqref{eq106} 
of equilibria at $x=\iota=1$ then imply that $h_1(k)$ is the $k$-th equilibrium, enumerated 
left to right, along the horizontal axis.
Suppose $h_1(k)=v_j$. Then, $h_1(k)=v_j=h_0(j)$, i.e.~the positional index $k$ of $v_j$ along the 
horizontal axis is given by the inverse Sturm permutation
\begin{equation} \label{eq200}
k=\sigma^{-1}(j)\,;
\end{equation}
see \eqref{eq107}. Consequently, in terms of the $h_1$-boundary order we obtain 
$v_j=h_1(\sigma^{-1}(j))$.

\begin{figure} 
\centering \includegraphics[scale = 1.0]{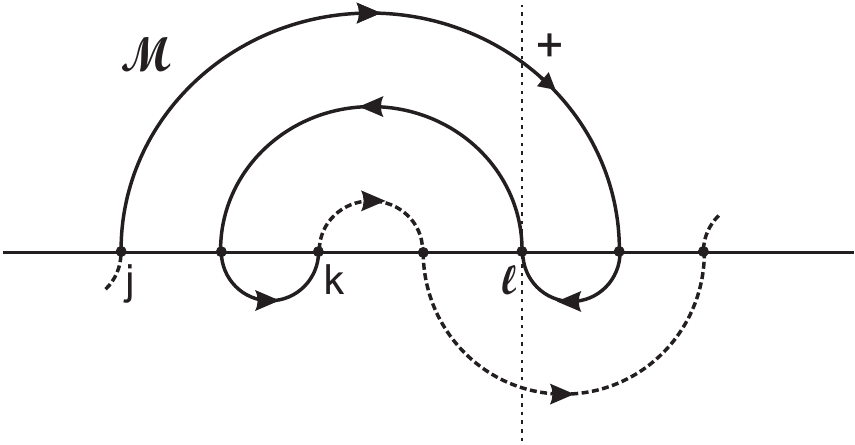}
\caption{\em\small The crossing number $c(j,k;\ell)$ counts the crossings between the meander 
segment of meander $\cM$ from $j$ to $k$ and the vertical line through $\ell$, clockwise with 
respect to $\ell$. Crossings at $\ell$ are not counted. Here $c(j,k;\ell)=1$. 
Note, $c(j,k;j)=c(j,k;k)=0$.}
\label{fig22}
\end{figure}

For any given Sturm meander $\cM$ in canonical form let $j<k$ and $\ell$ denote labels along 
$h_0$,~i.e. along the meander $\cM$, of three meander intersections with the horizontal axis; 
see Figure~\ref{fig22}.
Suppose the oriented meander segment $\cM_{j\,k}$ of $\cM$, from $j$ to $k$, intersects the vertical line 
through $\ell$ transversely. (We ignore any crossing at $\ell$ itself, in case $j\leq \ell\leq k$.)
We count clockwise crossings with respect to $\ell$ as $+1$, and counter-clockwise crossings as $-1$.
Then we define the {\em crossing number} 
$$ c=c(j,k;\ell) $$
of $\cM_{j\,k}$  with respect to $\ell$ as the total number of signed crossings, along the meander 
segment of $\cM$ from $j$ to $k$. 
In other words, $c(j,k;\ell)$ counts the total number of net clockwise half windings of the meander 
$\cM$ from $j$ to $k$, around $\ell$.
For $j=k$ we naturally define $c(j,j;\ell):=0$.

For nontransverse crossings with the vertical line through $\ell$,~e.g. during homotopies of meanders not 
in canonical form, the definition can be extended by transverse approximations. 
For details see \cite{firo99}, and our comment in the Appendix section \ref{sec5}.

For $j>k$, i.e.~if we follow the meander $\cM$ in reverse orientation, we define 
\begin{equation}\label{eq203} 
c(j,k;\ell) := -c(k,j;\ell)\,.
\end{equation}

Our definition of the half winding numbers $c(j,k;\ell)$ with respect to $\ell$ implies 
the additivity property
\begin{equation}\label{eq202} 
c(j_1,j_2;\ell) + c(j_2,j_3;\ell) = c(j_1,j_3;\ell)\,,
\end{equation}
for all choices of $j_1, j_2, j_3$, and all $\ell$.

The zero number $z(v_k-v_j)$ of the difference between any pair of distinct equilibria of any realization 
\eqref{eq102} of the Sturm meander $\cM$ can be obtained directly from the crossing numbers 
$c(j,k;\ell)$. Indeed, from \cite[Proposition~3]{roc91} it follows that
\begin{equation} \label{eq204}
z(v_k-v_j) = 
\begin{cases}
i(v_j) + c(j,k;j), & \text{ if } Q\cM_{j\,j+1} \text{ is odd} ; \\  
\\
i(v_j) - 1 + c(j,k;j), & \text{ if } Q\cM_{j\,j+1} \text{ is even} . 
\end{cases}
\end{equation}
Here $Q\cM_{j\,j+1}$ denotes the quadrant with respect to $j$ of the arc segment of $\cM$ strictly 
between the intersection points $j$ and $j+1$.
This is a nonlinear extension of the Sturm-Liouville property of solutions, which we review and prove
in the Appendix section~\ref{sec5}.

From \eqref{eq204} we recursively obtain the zero numbers $z_{j\,k}:=z(v_k-v_j)$ in terms of the meander 
permutation $\sigma$. In fact, additivity \eqref{eq202} implies 
\begin{equation} \label{eq207}
z_{j\,k+1} - z_{j\,k} = c(j,k+1;j) - c(j,k;j) = c(k,k+1;j) ,
\end{equation}
for all $1\le j<k\le N$. In terms of the permutation $\sigma$ we have
\begin{equation} \label{eq207a}
c(k,k+1;j) = \frac12(-1)^{k+1}\left[\sign\left(\sigma^{-1}(k+1)-\sigma^{-1}(j)\right) 
-\sign\left(\sigma^{-1}(k)-\sigma^{-1}(j)\right)\right] .
\end{equation}
Moreover, in view of dissipativeness of $f$ and the parabolic comparison principle, we obtain 
$v_1(x) < v_j(x) < v_N(x)$ for all $1<j<N$ and all $0\le x\le 1$. 
Hence, we define $z_{1\,j}=z_{j\,N}:=0$ for all $1<j<N$ and the zero numbers $z_{j\,k}$ for all 
$1\le j<k\le N$ are obtained from the decreasing recursion (see for example \cite{firo96})
\begin{equation} \label{eq207b}
\begin{aligned}
& z_{1\,j} = z_{j\,N} := 0 \\
& z_{j\,k} = z_{j\,k+1} + \frac12(-1)^{k}\left[\sign\left(\sigma^{-1}(k+1)-\sigma^{-1}(j)\right) 
-\sign\left(\sigma^{-1}(k)-\sigma^{-1}(j)\right)\right] .
\end{aligned}
\end{equation}

We emphasize that Morse indices $i(v_j)$ and zero numbers $z_{j\,k}$ provide all the information 
necessary to establish the existence of heteroclinic orbit connections between pairs of equilibria 
$v_j\leadsto v_k$. Such results derive from the zero number dropping argument mentioned in the 
Introduction, section \ref{sec1}, in relation with \eqref{eq108}; see \cite{firo96}. 
Here we recall the following notion of equilibria adjacency first introduced in \cite{wol02}. 
Two equilibria $v_j$ and $v_k$ are said to be {\em $z$-adjacent} if there does not exist any other 
equilibrium $w$ with $w(0)$ strictly between $v_j(0)$ and $v_k(0)$ such that 
\begin{equation} \label{eq208}
z(w-v_j)=z(v_k-w)=z(v_k-v_j) \,.
\end{equation} 
Then $v_j$ and $v_k$ are heteroclinically connected, $v_j\leadsto v_k$, if and only if $i(v_j)>i(v_k)$ and 
$v_j$ and $v_k$ are $z$-adjacent; see \cite[Theorem 2.1]{wol02} and also \cite[Appendix 7]{firo18b}.
If, due to an equilibrium $w$, the equilibria $v_j$ and $v_k$ are not $z$-adjacent, we say that $w$ 
satisfying \eqref{eq208} {\em blocks} heteroclinic connections between $v_j$ and $v_k$.

We now turn to the proof of our claims \eqref{eq116}--\eqref{eq116c}. 

\begin{lemma} \label{lem20}
Let $w=w^\iota_\pm$ denote a boundary neighbor of $\cO$, defined in \eqref{eq111}, satisfying $i(w)=n-1$. 
Then, $w\in\cE^{n-1}_\pm(\cO)$. More specifically, $w$ is the minimax equilibrium identified by 
\eqref{eq116}--\eqref{eq116c}.  
\end {lemma}

{\em Proof:} 
Indeed, if $\iota=0$, \eqref{eq204} implies $z(w-\cO)=n-1$. The same result also holds for the case 
of $\iota=1$ by invoking the involution $\tau: x \mapsto 1-x$. Moreover $\cO \leadsto w$, 
because such heteroclinic orbits are not blocked: see \eqref{eq208} and \cite{wol02}. 
This proves our claim $w\in \cE^{n-1}_\pm(\cO)$. 
Since $w$ is a meander neighbor of $\cO$, keeping track of $\sign(w^\iota_\pm-\cO)$ at the appropriate 
boundary $x=\iota=0,1$ for even/odd $n$, we obtain \eqref{eq116}--\eqref{eq116c}.
\hfill $\square$

\medskip

Finally we recall a result appearing in \cite[Lemma 1]{firowo12} (see also \cite{firo18b}) which relates 
the parity of the Morse indices $i(v_j)$ and the direction of the vertical crossing of the horizontal axis 
at $j$ by the meander $\cM$. 
Since it is used several times in our proof we highlight this result in the following lemma. 

\begin{lemma} \label{lem21}
Let $j=h^{-1}_0(v_j)$ denote the $h_0$-label of the $j$-th equilibrium $v_j\in\cE$ along the 
meander $\cM$ in canonical form. Then, the labels $j\in\{1,\dots,N\}$ and the Morse indices $i(v_j)$ have 
the opposite even/odd parity.
Similarly, if the Morse index $i(v_j)$ is even/odd then the meander $\cM$ crosses the horizontal 
line at $j$ vertically in the upwards/downwards direction, respectively.
\end{lemma}

{\em Proof:} For $j=1$ the first Morse index $i(v_1)=0$ is even and the first crossing of the 
meander $\cM$ with the horizontal axis at $1$ is vertical in the upwards direction. 
Afterwords, recursion \eqref{eq201} asserts that the even/odd parity of the Morse index 
$i(v_j)=i_j$ at the crossings $j$ alternates, along the meander $\cM$. Therefore, simple induction 
on $j$ along $\cM$ proves the lemma.
Similarly, the Jordan curve $\cM$ crosses the horizontal axis in the upwards and downwards 
direction, alternately. 
\hfill $\square$

\section{The minimax property}\label{sec3}

We now prove our main result Theorem~\ref{theo11}. Let again $\cO\in\cE$ denote the reference unstable 
equilibrium with $i(\cO)=n\ge 1$, and let $w^\iota_\pm, \iota\in\{0,1\}$, denote its four 
$h_\iota$-boundary neighbors given by \eqref{eq111}. 

In the Appendix section~\ref{sec6} we present {\em suspensions} of the Sturm meander $\cM$. 
Any suspension increases all the Morse indices of the equilibria by one, changes their even/odd parity, 
and adds two new equilibria to $\cE$. 
This operation corresponds to a double cone suspension of the global attractor $\cA$ by the 
addition of the new equilibria, but preserves the geometric structure of the previous global attractor. 
In particular, by Lemma~\ref{l62}, the minimax property is preserved by this suspension.  
Therefore, without loss of generality, in the following we assume $n$ to be odd. 

To outline our proof we remark that \eqref{eq117} in Theorem~\ref{theo11} consists of four different cases 
according to the four choices of the sign $\pm$ and the $\iota\in{0,1}$. We first consider the case of 
$\underline v^{\,1}_{\,+}$,~i.e. we choose the sign $+$ and consider the boundary $x=\iota=1$.

We assume that $i(w_+^1)=n-1$. Since $n$ is assumed odd, according to \eqref{eq116}--\eqref{eq116c} we 
have $\underline v^{\,1}_{\,+}=w^1_+\in \cE^{n-1}_+(\cO)$. 
Moreover, $i(\underline v^{\,1}_{\,+})=i(w_+^1)=n-1$ by assumption, and $n-1$ is even. 
Then it remains to show that
\begin{equation} \label{eq301}
\underline v^{\,1}_{\,+} = \overline v^{\,0}_{\,+} \,.
\end{equation} 
To prove this we use the $z$-adjacency notion \eqref{eq208} to show that the equilibrium 
$\underline v^{\,1}_{\,+}$ blocks heteroclinic connections from $\cO$ to any equilibrium $w\in\cE$ with 
$z(w-\cO)=(n-1)_+$ further away from $\cO$, at $x=0$, than $\underline v^{\,1}_{\,+}$. 
By \eqref{eq115} this will prove \eqref{eq301}.
We conclude the proof of Theorem~\ref{theo11} for the remaining three cases of \eqref{eq116} 
by applications of the trivial attractor equivalences generated by the involutions 
\begin{equation} \label{eq302}
\tau : x \mapsto 1-x \, , \quad \kappa : u \mapsto -u \,.
\end{equation}  

\vspace*{.3cm}

{\em Proof of Theorem~\ref{theo11} in the case: $\underline v^{\,1}_{\,+}=w_+^1$, $n$ odd.}

Fix an equilibrium $\cO$ such that $n:=i(\cO)\ge 1$ is odd. Let $j_0:=h^{-1}_0(\cO)$ denote the 
$h_0$-label of the equilibrium $\cO\in\cE$ along the meander $\cM$. Consequently, the $h_0$-label 
of the $x=0$ boundary successor $w^0_+$ of $\cO$,~i.e. of $w^0_+:=h_0(h^{-1}_0(\cO)+1)$, is given 
by $j_0+1=h^{-1}_0(w^0_+)=h^{-1}_0(\cO)+1$.
Also, we let $j_1:=h^{-1}_0(w^1_+)$ denote the $h_0$-label of the $x=1$ boundary successor $w^1_+$ 
of $\cO$,~i.e. of $w^1_+:=h_1(h^{-1}_1(\cO)+1)$. 
Hence $j_1=h^{-1}_0(w^1_+)=\sigma(h^{-1}_1(\cO)+1)$. 
By definition, the intersection points labeled $j_0$ and $j_1$ are $h_1$-neighbors along the 
horizontal axis, in this order.
Since $n$ is odd, the meander $\cM$ is oriented downwards at the intersection point corresponding 
to $\cO$, by Lemma~\ref{lem21}. Similarly, $\cM$ is oriented upwards at the intersection 
point corresponding to $w_+^1$.

\begin{figure}[t] 
\centering \includegraphics[scale = 0.85]{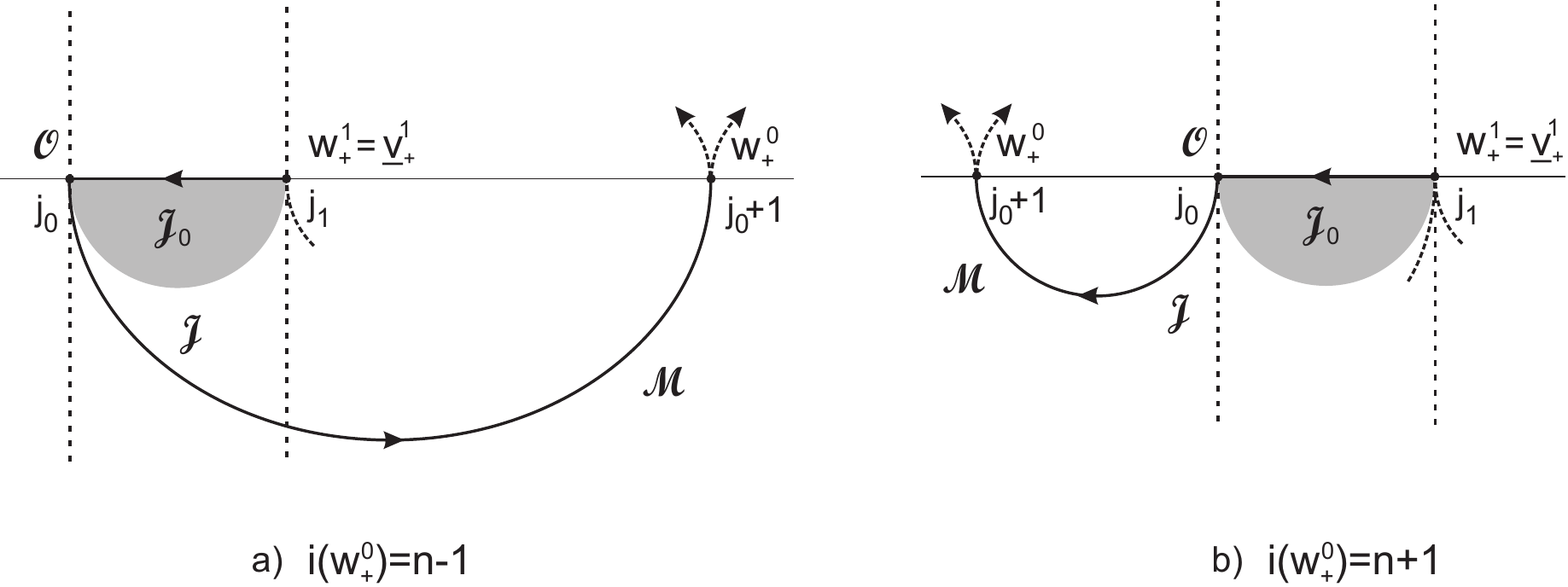}
\caption{\em\small Illustrations of the arc segment $\cM_{j_0\,j_0+1}$  of a meander $\cM$ 
with $i(\cO)=n$ odd and $\underline v^{\,1}_{\,+}=w^1_+\in\cE^{n-1}_+(\cO)$  under the alternative 
assumptions: a) $i(w_+^0)=n-1$, and b) $i(w_+^0)=n+1$.  
The union of the meander segment $\cM_{j_0\,j_1}$, from $j_0$ to $j_1$, with the section of the 
horizontal axis from $j_1$ back to $j_0$ is a planar Jordan curve $\cJ$.
The dark open semicircular disk $\cJ_0$ between the points labeled $j_0$ and $j_1$ is in the 
interior of the Jordan curve $\cJ$. }
\label{fig31}
\end{figure}

Our assumption $i(w_+^1)=n-1$ implies $\underline v^{\,1}_{\,+}=w^1_+\in\cE^{n-1}_+(\cO)$, 
by \eqref{eq114}, \eqref{eq116c}. Therefore \eqref{eq113} implies
\begin{equation} \label{eq303}
z_{j_0\,j_1} := z(w^1_+-\cO) = z(\underline v^{\,1}_{\,+}-\cO) = n-1 \,,
\end{equation} 
which is even. Hence, $\underline v^{\,1}_{\,+}(x)=w^1_+(x)>\cO(x)$ at both boundaries 
$x=\iota\in\{0,1\}$.
See Figure~\ref{fig31} for the illustration of these notations for meanders $\cM$ with 
$\underline v^{\,1}_{\,+}=w^1_+\in\cE^{n-1}_+(\cO)$ and odd $n$ in the alternative cases 
$i(w_+^0)=n\pm 1$ arising from $i(\cO)=n$ and \eqref{eq111}.

Let $\cJ$ denote the union of the meander segment $\cM_{j_0\,j_1}$, oriented from $j_0$ to $j_1$, with the 
short segment of the horizontal axis from $j_1$ back to the adjacent intersection $j_0$. Since $\cM$ is a 
meander and $j_0, j_1$ are horizontally adjacent, the closed oriented curve $\cJ$ is a planar Jordan curve. 
Moreover, as the Morse indices count the clockwise half-windings of the unit tangent vector along the 
canonical meander $\cM$ (see comment to \eqref{eq201}), the assumption 
$i(\underline v^{\,1}_{\,+})=i(w^1_+)=i(\cO)-1$ implies that $\cJ$ performs 
a full counter-clockwise winding, that is $\cJ$ is positively (i.e. left, counterclockwise) oriented.
This shows that the (dark) open semicircular disk $\cJ_0$ in the lower half-plane, with diameter given by 
the horizontal axis from $j_0$ to $j_1$, is interior to $\cJ$. See again Figure~\ref{fig31}.
In addition, the unbounded continuing segment of the meander $\cM$ which starts at $j_1$ and runs 
to the Northeast (eventually) asymptotically (see the meander definition in the beginning of section 
\ref{sec2}), is exterior to $\cJ$.

We now prepare to show claim \eqref{eq301} by contradiction. The strategy here is to suppose, 
contrary to \eqref{eq301}, that 
\begin{equation} \label{eq305}
\overline v^{\,0}_{\,+} \ne \underline v^{\,1}_{\,+} \,.
\end{equation}
In step 1 we then prove 
\begin{equation} \label{eq306}
z_{j_0\,j_2} := z(\overline v^{\,0}_{\,+}-\cO) = n-1 \,.
\end{equation}
Step 2 will establish that 
\begin{equation} \label{eq307}
z_{j_1\,j_2} := z(\overline v^{\,0}_{\,+}-\underline v^{\,1}_{\,+}) = n-1 \,.
\end{equation}
In view of \eqref{eq303} this will imply that equilibria $\cO$ and $\overline v^{\,0}_{\,+}$ cannot be 
$z$-adjacent, because $\underline v^{\,1}_{\,+}$ will block $\cO\leadsto\overline v^{\,0}_{\,+}$; 
see \eqref{eq208}. This provides a contradiction to definition \eqref{eq113}, \eqref{eq115} of 
$\overline v^{\,0}_{\,+}$. 
These three steps will prove claim \eqref{eq301}, for odd $n=i(\cO)$ and $\underline v^{\,1}_{\,+}=w^1_+$ 
with $i(w^1_+)=n-1$.

Let us initiate this indirect strategy, therefore, by assuming that the $h_0$-maximal element 
$\overline v^{\,0}_{\,+}$ of $\cE^{n-1}_+(\cO)$, most distant from $\cO$ at $x=0$, differs from 
the $h_1$-minimal element $\underline v^{\,1}_{\,+}$ in $\cE^{n-1}_+(\cO)$, closest to $\cO$ at $x=1$, 
i.e. \eqref{eq305}.

\vspace*{.3cm}

{\em Step 1: Proof of \eqref{eq306}}

\begin{figure}[t] 
\centering \includegraphics[scale = 0.85]{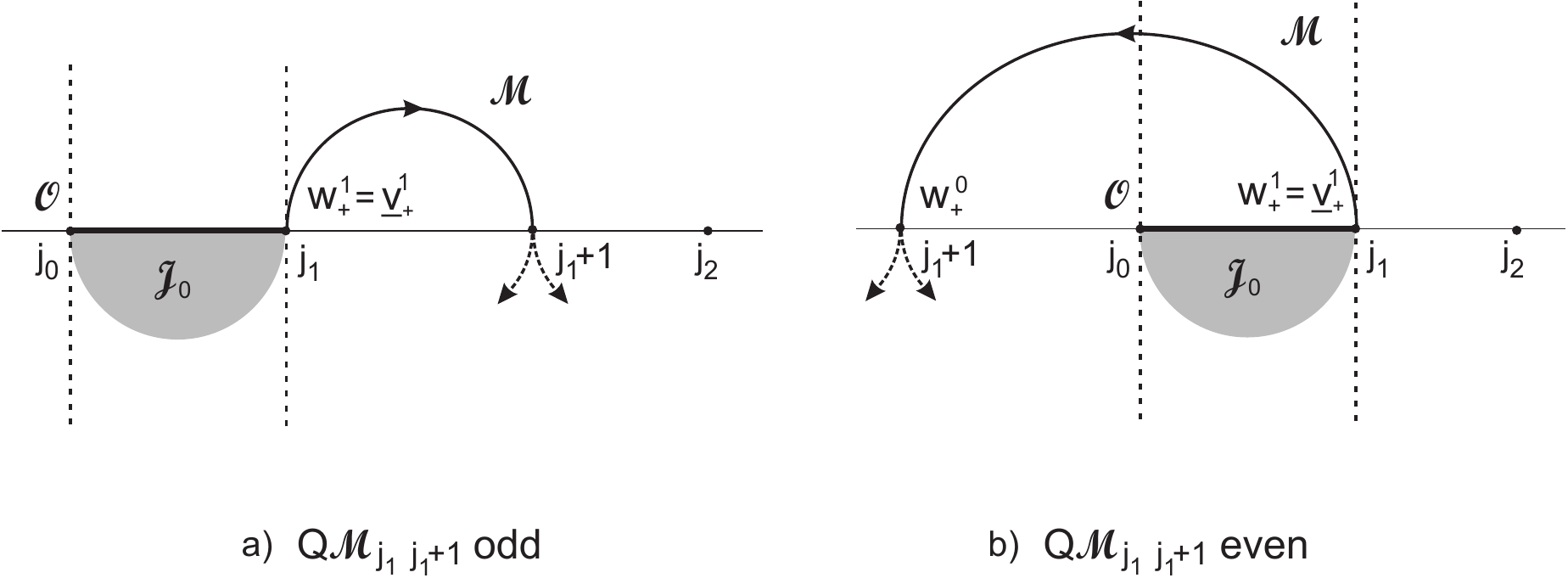}
\caption{\em\small Illustrations of arc segments $\cM_{j_1\,j_1+1}$ of a meander $\cM$ 
with $i(\cO)=n$ odd and $\underline v^{\,1}_{\,+}=w^1_+\in\cE^{n-1}_+(\cO)$ under the alternative 
assumptions: a) $Q\cM_{j_1\,j_1+1}$ is odd, and b) $Q\cM_{j_1\,j_1+1}$ is even. 
The meander segment $\cM_{j_1\,j_2}$ starting at the point $j_1$ is in the exterior of $\cJ$ and 
cannot cross the horizontal axis between $j_0$ and $j_1$. In case a) we obtain 
$c(j_1,j_1+1;j_1)=c(j_1,j_1+1;j_0)$, and in case b) we obtain $c(j_1,j_1+1;j_1)=c(j_1,j_1+1;j_0)+1$.}
\label{fig32}
\end{figure}

Let $j_2$ denote the $h_0$-label of the $h_0$-maximal equilibrium 
$\overline v^{\,0}_{\,+}\in\cE^{n-1}_+(\cO)$ most distant from $\cO$ at $x=0$,~i.e. 
$v_{j_2}:=\overline v^{\,0}_{\,+}=h_0(j_2)$. 
Since $\overline v^{\,0}_{\,+}$ is $h_0$-maximal in $\cE^{n-1}_+(\cO)$, by \eqref{eq114}, and we 
just assumed $\overline v^{\,0}_{\,+}\ne \underline v^{\,1}_{\,+}=w^1_+\in\cE^{n-1}_+(\cO)$, we 
have $j_2>j_1$. Definition \eqref{eq113} of $\cE^{n-1}_+(\cO)$ implies \eqref{eq306} as claimed by our 
strategy. Moreover, since $n$ is assumed odd, $n-1$ is even and $v_{j_2}>\cO$ at the boundary $x=1$. 
See Figure~\ref{fig32}.

{\em Step 2: Proof of \eqref{eq307}}

We first recall that \eqref{eq303} and \eqref{eq306} assert $z_{j_0\,j_1}=n-1=z_{j_0\,j_2}$. 
Invoking \eqref{eq204} with $j=j_0$ and $k=j_1,j_2$ we obtain equality of the crossing numbers
\begin{equation} \label{eq308}
c(j_0,j_1;j_0) = c(j_0,j_2;j_0) \,.
\end{equation}
Indeed $j_1$ and $j_2$ both occur after $j_0$ on the meander $\cM$, and therefore refer to the same 
even/odd quadrant of $\cM_{j_0\,j_0+1}$.
Since $c(j_0,j_2;j_0)=c(j_0,j_1;j_0)+c(j_1,j_2;j_0)$, by the additivity property \eqref{eq202}, 
this shows 
\begin{equation} \label{eq309}
c(j_1,j_2;j_0) = 0 \,.
\end{equation}

Next, since the intersection points $j_0$ and $j_1$ are $h_1$-adjacent neighbors along the 
horizontal axis, any semicircular arc segment of the meander segment $\cM_{j_1+1\,j_2}$ which 
crosses the vertical line at $j_1$ also crosses the vertical line of $j_0$, in the same upper 
or lower half-plane and in the same direction. 
The total contribution of $\cM_{j_1+1\,j_2}$ to those crossing numbers therefore satisfies
\begin{equation} \label{eq310}
c(j_1+1,j_2;j_0) = c(j_1+1,j_2;j_1) \,.
\end{equation}

For the first segment $\cM_{j_1\,j_1+1}$ we consider the alternative cases of even/odd $Q\cM_{j_1\,j_1+1}$ 
and obtain
\begin{equation} \label{eq311}
c(j_1,j_1+1;j_1) = 
\begin{cases}
c(j_1,j_1+1;j_0), & \text{ if } Q\cM_{j_1\,j_1+1} \text{ is odd} \,; \\
\\
c(j_1,j_1+1;j_0) + 1, & \text{ if } Q\cM_{j_1\,j_1+1} \text{ is even} \,.
\end{cases}
\end{equation}
Indeed, the crossing number definition implies (see Figure~\ref{fig32}):
\begin{equation} \label{eq311a}
c(j_1,j_1+1;j_1) = 0 \,, \quad c(j_1,j_1+1;j_0) = 
\begin{cases}
0, & \text{ if } Q\cM_{j_1\,j_1+1} \text{ is odd} \,; \\
\\
-1, & \text{ if } Q\cM_{j_1\,j_1+1} \text{ is even} \,.
\end{cases}
\end{equation}

From \eqref{eq310}, \eqref{eq311}, by additivity \eqref{eq202} and using \eqref{eq309} we have 
\begin{equation} \label{eq312}
c(j_1,j_2;j_1) = 
\begin{cases}
c(j_1,j_2;j_0) = 0, & \text{ if } Q\cM_{j_1\,j_1+1} \text{ is odd} \,; \\
\\
c(j_1,j_2;j_0) + 1 = 1, & \text{ if } Q\cM_{j_1\,j_1+1} \text{ is even} \,.
\end{cases}
\end{equation}
Therefore, by \eqref{eq204}, and since $i(\underline v^{\,1}_{\,+})=i(w^1_+)=n-1$ by assumption, 
from \eqref{eq312} we conclude
\begin{equation} \label{eq313}
z_{j_1\,j_2} = 
\begin{cases} 
i(\underline v^{\,1}_{\,+})+c(j_1,j_2;j_1) = n-1, & \text{ if } Q\cM_{j_1\,j_1+1} \text{ is odd} \,; \\
\\
i(\underline v^{\,1}_{\,+})-1+c(j_1,j_2;j_1) = n-1, & \text{ if } Q\cM_{j_1\,j_1+1} \text{ is even} \,.
\end{cases}
\end{equation}
This proves claim \eqref{eq307}.

\vspace*{.3cm}

{\em Step 3: Proof of blocking}

To reach a contradiction to our indirect assumption \eqref{eq305}, we now show that the equilibrium 
$w:=\underline v^{\,1}_{\,+}$ blocks the heteroclinic connection 
$\cO\leadsto\overline v^{\,0}_{\,+}\in\cE^{n-1}_+(\cO)$, by blocking property \eqref{eq208} and contrary 
to definition \eqref{eq113} of $\cE^{n-1}_+(\cO)$. 
Indeed, relations \eqref{eq303} for $w^1_+=\underline v^{\,1}_{\,+}$, \eqref{eq306}, and \eqref{eq307} assert 
\begin{equation} \label{eq315}
z(\underline v^{\,1}_{\,+}-\cO) = n-1 \,, \quad z(\overline v^{\,0}_{\,+}-\cO) = n-1 \,, \quad
z(\overline v^{\,0}_{\,+}-\underline v^{\,1}_{\,+}) = n-1 \,,
\end{equation}
respectively. Moreover, the ordering $j_0<j_1<j_2$ corresponds to the ordering of the equilibria 
$\cO=h_0(j_0)$, $\underline v^{\,1}_{\,+}=h_0(j_1)$, $\overline v^{\,0}_{\,+}=h_0(j_2)$ at the $h_0$-boundary:
\begin{equation} \label{eq316}
\cO < \underline v^{\,1}_{\,+} < \overline v^{\,0}_{\,+} \quad \mbox{ at } x=0 \,.
\end{equation}
Therefore, \eqref{eq208} and \eqref{eq315} show that the equilibrium $\underline v^{\,1}_{\,+}$ blocks 
heteroclinic connections $\cO\leadsto\overline v^{\,0}_{\,+}$. This contradiction to the definition of 
$\overline v^{\,0}_{\,+}$ in \eqref{eq115}, \eqref{eq113}, proves claim \eqref{eq301} for 
$\underline v^{\,1}_{\,+}=w_+^1$. 
\hfill $\square$

\vspace*{.3cm}

{\em Proof of Theorem~\ref{theo11} in the remaining three cases.}

By the trivial equivalence $\tau: x\mapsto 1-x$ we obtain a Sturm global attractor $\tau\cA$ with the 
equilibria $\tau\underline v^{\,\iota}_{\,+}, \tau\overline v^{\,\iota}_{\,+}, \iota\in\{0,1\}$, now 
referring to $\tau\cO$. Specifically,
\begin{equation}\label{eq317} 
\tau\underline v^{\,1-\iota}_{\,+} = \underline v^{\,\iota}_{\,+} \ ; \quad 
\tau\overline v^{\,1-\iota}_{\,+} = \overline v^{\,\iota}_{\,+} \ ; \quad \iota\in\{0,1\} \,.
\end{equation}

Then, if the appropriate $i(\tau w_+^\iota)=n-1$, by \eqref{eq301} we have that 
$\tau\underline v^{\,1}_{\,+}=\tau\overline v^{\,0}_{\,+}$. This shows that, \eqref{eq317} implies
\begin{equation}\label{eq318} 
\underline v^{\,0}_{\,+} = \overline v^{\,1}_{\,+} \,,
\end{equation}
which settles the remaining case with sign $+$.

For the equilibria in $\cE^{n-1}_-(\cO)$, if we assume that $i(w_-^\iota)=n-1$ for the corresponding 
boundary predecessors $w_-^\iota$ with $\iota\in\{0,1\}$, the remaining two 
cases for the sign $-$ follow by the trivial equivalence $\kappa: u\mapsto -u$. 

For the convenience of the reader we include a table of the action of the Klein group $\Z_2\otimes\Z_2$ 
generated by the commuting involutions $\langle\tau,\kappa\rangle$ on the equilibria 
$\underline v^{\,\iota}_{\,\pm}$.  

\begin{center}
\begin{tabular}{ |c|c|c|c| } 
 \hline
 & & & \\
 $\cO$ & $\tau\cO$ & $\kappa\cO$ & $\tau\kappa\cO$ \\ 
 & & & \\
 \hline
 & & & \\
 $\underline v^{\,\iota}_{\,\pm}(\cO)$ & $\underline v^{\,1-\iota}_{\,\pm}(\tau\cO)$ &
 $\underline v^{\,\iota}_{\,\mp}(\kappa\cO)$ & $\underline v^{\,1-\iota}_{\,\mp}(\tau\kappa\cO)$ \\ 
 & & & \\
 \hline
\end{tabular}
\end{center}

This completes the proof of our main Theorem~\ref{theo11}. 
\hfill $\square$

\section{Discussion and Example} \label{sec4}

As pointed out in the Introduction, section \ref{sec1}, the minimax property \eqref{eq117} 
simplifies the task of identifying the equilibria $\cE^{n-1}_\pm(\cO)$ of \eqref{eq113} directly 
from the Sturm meander $\cM$.

\begin{figure}[t] 
\centering \includegraphics[scale = 0.8]{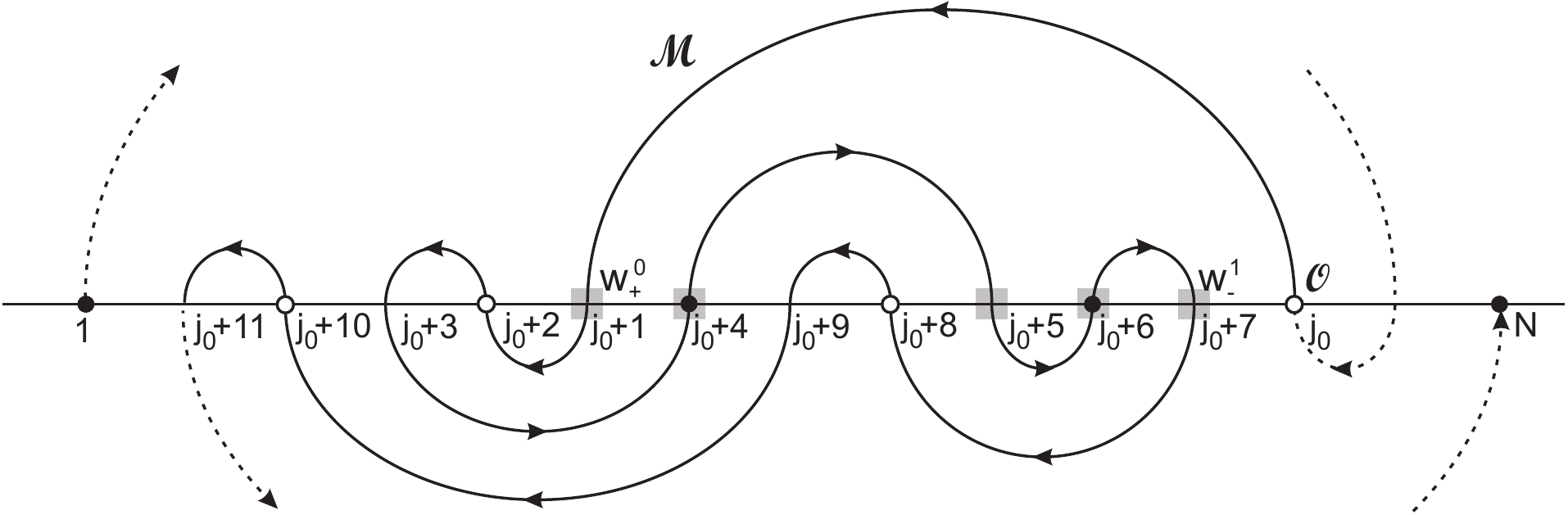}
\caption{\em\small Illustration of the meander section $\cM_{j_0\,j_0+11}$ corresponding to the 
partial Sturm permutation $\sigma = \{ 1 \ \dots \ j_0+11 \ j_0+10 \ j_0+3 \ j_0+2 \ j_0+1 \ 
j_0+4 \ j_0+9 \ j_0+8 \ j_0+5 \ j_0+6 \ j_0+7 \ j_0 \ \dots \ N \}$, 
where the reference unstable equilibrium $\cO=h_0(j_0)$ has even Morse index $n:=i(\cO)=2$ and odd 
$i(\underline v^{\,1}_{\,+})=i(w^1_-)=i(w^0_+)=i(\cO)-1=n-1=1$. 
The equilibria with Morse index $i=0$ are indicated by black dots and the equilibria with Morse 
index $i=2$ are indicated by white dots. 
Equilibria with Morse index $i=1$ are indicated by simple intersections. The five intersections at 
$j_0+\{1, 4, 5, 6, 7\}$ in gray squares indicate the target set $\cE^{n-1}_+(\cO)$ defined in 
\eqref{eq113}.} 
\label{fig41}
\end{figure}

To emphasize the ``local'' aspect of our global result, we next show an example for which the 
identification of the equilibria in $\cE^{i(\cO)-1}_+(\cO)$ is obtained from the knowledge of only 
a segment of the meander $\cM$. In fact, we will only prescribe a segment of the Sturm permutation 
$\sigma$.
We assume that our reference equilibrium $\cO\in\cE$ has even unstable dimension $n:=i(\cO)=2$. As 
before, $j_0=h^{-1}_0(\cO)$ denotes the $h_0$-label of $\cO$ along the meander. 
We consider a Sturm permutation according to the following template
\begin{equation}\label{eq401} 
\sigma = \{ 1 \ \dots \ j_0+11 \ j_0+10 \ j_0+3 \ j_0+2 \ j_0+1 \ j_0+4 \ j_0+9 \ j_0+8 \ j_0+5 \ 
j_0+6 \ j_0+7 \ j_0 \ \dots \ N \} \,.
\end{equation}
This corresponds to a meander section $\cM_{j_0\,j_0+11}$ which we illustrate in 
Figure~\ref{fig41}. 
Note the orientation of the meander $\cM$ due to the assumption of even $n$ and Lemma~\ref{lem21}. 
Our objective is to identify the set of equilibria $\cE^{n-1}_+(\cO)=\cE^{1}_+(\cO)$ from the 
partial ``local'' information \eqref{eq401} on $\sigma$. 

By \eqref{eq111} the boundary neighbors $w^0_+=w^0_+(\cO)$ and $w^1_-=w^1_-(\cO)$ are given by
\begin{equation}\label{eq402} 
w^0_+ = h_0(h^{-1}_0(\cO)+1) = h_0(j_0+1) \ ; \qquad w^1_- = h_1(h^{-1}_1(\cO)-1) = h_0(j_0+7) \,.
\end{equation}
Indeed, $w^1_- = h_0(j_0+7)$ since $w^1_-$ is the left $h_1$-neighbor of $\cO$.
The recursion \eqref{eq201} for the Morse indices $i_j$ applied to $\sigma$ in \eqref{eq401}, implies 
$i(w^1_-)=i(w^0_+)=i(\cO)-1=n-1=1$. 
Therefore, Theorem~\ref{theo11}, \eqref{eq117}, together with \eqref{eq116a}, \eqref{eq116b}  yields 
\begin{equation}\label{eq403} 
w^0_+ = \underline v^{\,0}_{\,+} = \overline v^{\,1}_{\,+} \ ; \quad w^1_- = \underline v^{\,1}_{\,+} = 
\overline v^{\,0}_{\,+} \,.
\end{equation}
Equation \eqref{eq403} implies that the set $\cE^1_+(\cO)$ of \eqref{eq113} is strictly 
contained in the set of equilibria $v_j$ with values at the boundary $x=0$ in the interval
\begin{equation}\label{eq404} 
\cO(0) < w^0_+(0) = \underline v^{\,0}_{\,+}(0) = v_{j_0+1}(0) \le v_j(0) \le v_{j_0+7}(0) = 
\overline v^{\,0}_{\,+}(0) = w^1_-(0) \,,
\end{equation} 
that is, the set of all equilibria $v_j=h_0(j)$ with $j_0+1\le j\le j_0+7$. 
Equation \eqref{eq403} also implies that $\cE^1_+(\cO)$ is strictly contained in the set of 
equilibria with values at the boundary $x=1$ in the interval
\begin{equation}\label{eq405} 
w^0_+(1) = \overline v^{\,1}_{\,+}(1) = v_{j_0+1}(1) \le v_j(1) \le v_{j_0+7}(1) = 
\underline v^{\,1}_{\,+}(1) = w^1_-(1) < \cO(1) \,.
\end{equation}

We now claim that
\begin{equation}\label{eq406} 
\cE_+^1(\cO) =: \{v_{j_0+k}: k\in K\}, \quad \textrm{with} \quad K=\{1,4,5,6,7\}\,.
\end{equation}
Here the index set $K$ is defined by the left equality.

Equation \eqref{eq404} implies $K\subseteq \{1,\ldots,7\}$.
Equation \eqref{eq405} implies $K\subseteq \{1,4,9,8,5,6,7\}$.
By intersection, this proves $K\subseteq \{1,4,5,6,7\}$.
 
We have to show, conversely, that $K\supseteq \{1,4,5,6,7\}$.
The zero number formula \eqref{eq204}, with $j=j_0, \,v_j=\cO$, and $v_{j_0+k}$ replacing $v_k$ 
there, implies $z(v_{j_0+k}-\cO)=1_+$\,, for any $k\in \{1,4,5,6,7\}$. 
Therefore it only remains to prove 
\begin{equation}\label{eq406a} 
k\in \{1,4,5,6,7\} \quad \Longrightarrow \quad \cO \leadsto v_{j_0+k}\,.
\end{equation}

To prove claim \eqref{eq406a}, we recall that $v_j\leadsto v_k$ are heteroclinically connected, 
if and only if $i(v_j)>i(v_k)$ and $v_j$ and $v_k$ are $z$-adjacent; see \eqref{eq208}. 
The following Morse indices $i(v_j)=i_j$ are easily determined from \eqref{eq201}: 
we have $i(v_{j_0+k})=0$, for $k\in\{4,6\}$, and $i(v_{j_0+k})=1$, for $k\in\{1,5,7\}$.

Two equilibria $v_j$ and $v_{j+1}$,~i.e. $h_0$-boundary neighbors, are automatically $z$-adjacent. 
In fact, \eqref{eq201} implies $|i(v_j)-i(v_{j+1})|=1$ and there is no third equilibrium $w$ with 
$w(0)$ strictly between $v_j(0)$ and $v_{j+1}(0)$. 
In particular $\cO=v_{j_0}\leadsto v_{j_0+1}=w^0_+$\,. 
The same assertion holds for two equilibria $v_{\sigma^{-1}(j)}$ and $v_{\sigma^{-1}({j+1})}$,~i.e. 
$h_1$-boundary neighbors at $x=1$, by the trivial equivalence $\tau: x \mapsto 1-x$. 
In particular $\cO=v_{j_0}\leadsto v_{j_0+7}=w^1_-$\,. This takes care of the cases $k\in\{1,7\}$.

From \eqref{eq401} (see also the template Figure~\ref{fig41}), we immediately verify the $h_0$-boundary 
$z$\nobreakdash-adjacency of the equilibria $v_{j_0+7}, v_{j_0+6}$\,, and the $h_1$-boundary 
$z$\nobreakdash-adjacency of the equilibria $v_{j_0+1}, v_{j_0+4}$. 
Hence, we obtain $w^1_-=v_{j_0+7}\leadsto v_{j_0+6}$ and $w^0_+=v_{j_0+1}\leadsto v_{j_0+4}$. 
By heteroclinic transitivity, this takes care of claim \eqref{eq406a} for $k\in\{4,6\}$.

For the final equilibrium $v_{j_0+k}$ with $k=5$, we consider all potentially $z=1$ blocking 
equilibria $w$ with $w(1)$ strictly between $\cO(1)$ and $v_{j+5}(1)$. By template \eqref{eq401}, 
there are exactly two such candidate equilibria: $w=v_{j_0+6}$ and $w=v_{j_0+7}$. Now, using the 
bottom row of \eqref{eq204} with $j=j_0+5, k=j_0+6, j_0+7$, we obtain the zero numbers
\begin{equation}\label{eq407} 
z_{j_0+6\,j_0+5} = z(v_{j_0+5}-v_{j_0+6}) = 0\,, \quad z_{j_0+7\,j_0+5} = z(v_{j_0+5}-v_{j_0+7}) 
= 0 \,.
\end{equation}
This shows that neither $w=v_{j_0+6}$\,, nor $w=v_{j_0+7}$\,, can block the $z=1$ heteroclinic 
connection $\cO\leadsto v_{j_0+5}$, at $x=1$, and we conclude that also $v_{j_0+5}\in\cE^1_+(\cO)$ 
by $z$-adjacency.
This completes the proof of our claims \eqref{eq406a} and \eqref{eq406}.

\begin{figure}[t] 
\centering \includegraphics[scale = 1.0]{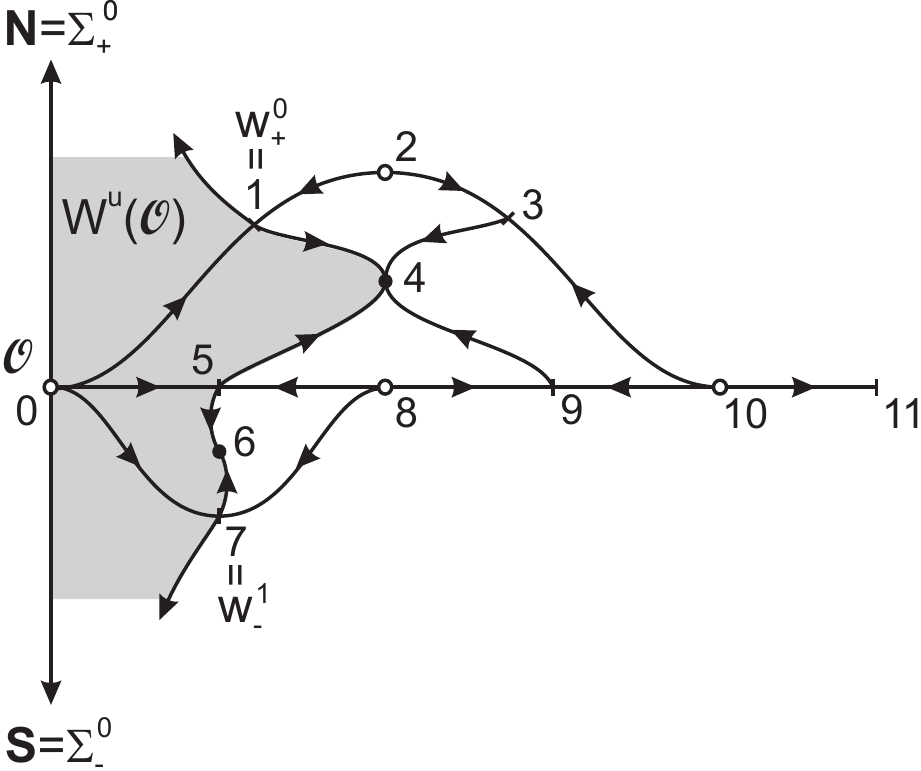}
\caption{\em\small Sketch of the heteroclinic connections in the global attractor corresponding to 
the partial Sturm permutation $\sigma$. To abbreviate the equilibrium labels we let $j_0=0$. 
The labels corresponding to $\cO$, $w^0_+$ and $w^1_-$ are $0$, $1$ and $7$, respectively. 
The shaded area corresponds to the subset of $W^u(\cO)$ composed of all heteroclic connections 
$\cO\leadsto w$ with $z(w-\cO)=2_+$. In $W^u(\cO)$ there are exactly two heteroclinic connections 
$\cO\leadsto w$ with $z(w-\cO)=0_\pm$ which select the unique polar equilibria 
${\bf N}=\Sigma^0_-(\cO)$ and ${\bf S}=\Sigma^0_+(\cO)$; see \cite{firo18b}. 
These polar equilibria are not labeled by the partial Sturm permutation $\sigma$. According to the 
dashed template in Figure~4.1, one possibility is ${\bf N}=-2$ and ${\bf S}=12$.} 
\label{fig42}
\end{figure}

In Figure~\ref{fig42} we sketch all the single-orbit heteroclinic connections obtained from the 
partial Sturm permutation $\sigma$. 
In particular, the set $\cE^1_+(\cO)$ is illustrated with its chain of heteroclinic connections  
between the boundary neighbors $w^0_+$ and $w^1_-$.

Figure~\ref{fig42} displays all the single-orbit heteroclinic connections as obtained from 
$\sigma$, by checking $z$-adjacency directly from the recurrence \eqref{eq207b}. For completion we 
include next the partial zero number matrix corresponding to the Sturm permutation segment $\sigma$, 
completed on the diagonal by the Morse indices of the equilibria:
\begin{equation}\label{eq408} 
\left[z_{j\,k}\right]_{j,k=j_0,\dots,j_0+11} = 
\left( 
\begin{array}{cccccccccccc}
2 & 1 & 1 & 1 & 1 & 1 & 1 & 1 & 1 & 1 & 1 & 1 \\
1 & 1 & 1 & 1 & 0 & 0 & 0 & 0 & 0 & 0 & 1 & 1 \\
1 & 1 & 2 & 1 & 0 & 0 & 0 & 0 & 0 & 0 & 1 & 1 \\
1 & 1 & 1 & 1 & 0 & 0 & 0 & 0 & 0 & 0 & 1 & 1 \\
1 & 0 & 0 & 0 & 0 & 0 & 0 & 0 & 0 & 0 & 1 & 1 \\
1 & 0 & 0 & 0 & 0 & 1 & 0 & 0 & 1 & 1 & 1 & 1 \\
1 & 0 & 0 & 0 & 0 & 0 & 0 & 0 & 1 & 1 & 1 & 1 \\
1 & 0 & 0 & 0 & 0 & 0 & 0 & 1 & 1 & 1 & 1 & 1 \\
1 & 0 & 0 & 0 & 0 & 1 & 1 & 1 & 2 & 1 & 1 & 1 \\
1 & 0 & 0 & 0 & 0 & 1 & 1 & 1 & 1 & 1 & 1 & 1 \\
1 & 1 & 1 & 1 & 1 & 1 & 1 & 1 & 1 & 1 & 2 & 1 \\
1 & 1 & 1 & 1 & 1 & 1 & 1 & 1 & 1 & 1 & 1 & 1 
\end{array}
\right) \,.
\end{equation}

\section{Appendix: Nonlinear Sturm-Liouville property}\label{sec5}

In this Appendix we review and prove the nonlinear Sturm-Liouville property (NSL for short) in our meander 
setting. 
Let $v_j=h_0(j), v_k=h_0(k)$ denote two equilibria with $j<k$.
Then, the  NSL property corresponds to the relation \eqref{eq204} between zero numbers, Morse indices 
and crossing numbers, which we now repeat for convenience. 
The claim is that the zero number $z_{j\,k}:=z(v_k-v_j)$ is given by (\cite[Proposition~3]{roc91})
\begin{equation} \label{eq501}
z(v_k-v_j) = 
\begin{cases}
\ i(v_j) + c(j,k;j), & \text{ if } Q\cM_{j\,j+1} \text{ is odd } \,; \\  
\\
\ i(v_j) - 1 + c(j,k;j), & \text{ if } Q\cM_{j\,j+1} \text{ is even } \,, 
\end{cases}
\end{equation}
Here $Q\cM_{j\,j+1}$ denotes the quadrant with respect to $j$ of the arc segment of $\cM$ between the 
intersection points $j$ and $j+1$.

As in section \ref{sec2}, $c(j,k;j)$ denotes the net signed clockwise crossings of the oriented 
meander segment $\cM_{j\,k}$ from equilibrium crossing $j$ to $k$ through the vertical line of $j$, ignoring 
that first crossing. See \eqref{eq203}, \eqref{eq202}.

Our proof of the NSL property \eqref{eq501} is based on zeros and winding numbers associated to the 
solutions $v=v(x,a)$ of the initial value second order ODE problem
\begin{equation} \label{eq502}
0=v_{xx}+f(x,v,v_x) \quad , \quad v(0,a)=a , \quad v_x(0,a)=0 .
\end{equation}
The equilibrium boundary value problem \eqref{eq102} is related to this ODE by the shooting condition 
$v_x(x,a)=0$ at the right boundary $x=1$.
Let $a_j:=v_j(0),\ a_k:=v_k(0)$ denote the initial values at $x=0$ of the equilibria $v_j, v_k\in\cE$, i.e.
\begin{equation} \label{eq503}
v(\cdot,a_j)=v_j(\cdot) \,, \qquad v(\cdot,a_k)=v_k(\cdot) \,.
\end{equation}
Then the segment $\cM_{j\,k}$ of the meander $\cM$, from the intersection point $j$ to the 
intersection point $k$, is given by the planar curve
\begin{equation} \label{eq504}
a\mapsto (v(1,a),v_x(1,a))\in\R^2 \,, \qquad a_j \le a \le a_k 
\end{equation}
of boundary values at the shooting boundary $x=1$.
Note how the shooting condition $v_x(1,a)=0$ is actually satisfied, precisely, at the equilibrium 
intersections of the meander $\cM$ with the horizontal axis $v_x=0$ in the $(v,v_x)$-plane.

For $a_j<a\le a_k$\,, let 
\begin{equation}
\label{eq5w}
w=w(x,a):=(v(x,a)-v_j(x))/(a-a_j)
\end{equation}
denote the scaled difference between the two solutions $v,v_j$ of \eqref{eq502}.
Note that $w=w(x,a)$ solves a linear second order ODE initial value problem
\begin{equation} \label{eq507}
0 = w_{xx}+q_1(x,a)w_x+q_0(x,a)w \,, \qquad w(0,a)=1, \quad w_x(0,a)=0 \,.
\end{equation}
Indeed, the coefficients $q_0, q_1$ depend on $v, v_j$ and are given explicitly as
\begin{equation} \label{eq508}
\begin{aligned} 
q_0(x,a) = & \int_0^1 \partial_u f(x,r(x,a,\mu),r_x(x,a,\mu)) \ d\mu \,, \\
q_1(x,a) = & \int_0^1 \partial_p f(x,r(x,a,\mu),r_x(x,a,\mu)) \ d\mu \,,
\end{aligned}
\end{equation}
with $r(x,a,\mu) = \mu v(x,a) + (1-\mu) v_j(x)$, by the Fundamental Theorem of Calculus.

To extend the above construction \eqref{eq5w}--\eqref{eq508}, down to $a=a_j$\,, let
\begin{equation} \label{eq509}
q_0(x,a_j) := \partial_u f(x,v_j(x),v_{j\,x}(x)) \,, \qquad 
q_1(x,a_j) := \partial_p f(x,v_j(x),v_{j\,x}(x)) \,.
\end{equation}
Then $f\in C^1$ implies continuity of $q_0, q_1$\,, and therefore continuity of $w,w_x\in C^1$, 
in the closed rectangle 
\begin{equation}
\label{eq5R}
(x,a)\in\mathcal{R}:=[0,1]\times[a_j,a_k]\,.
\end{equation}

The initial condition $w(0,a)=1$ of the linear equation \eqref{eq507} implies $(w(x,a),w_x(x,a))\neq (0,0)$ 
on $\mathcal{R}$. We can therefore introduce (clockwise!) polar coordinates, according to the Pr\"ufer 
transformation
\begin{equation}
\label{eq5theta}
w=\rho\cos\vartheta\,, \quad w_x=-\rho\sin\vartheta\,
\end{equation}
and obtain
\begin{equation} \label{eq514}
\begin{aligned}
\rho_x =& - \rho \left( q_1(x,a) \sin^2\vartheta +(1-q_0(x,a)) \sin\vartheta\cos\vartheta \right) \\
\\
\vartheta_x =& \sin^2\vartheta - q_1(x,a) \sin\vartheta\cos\vartheta + q_0(x,a) \cos^2\vartheta \,.
\end{aligned}
\end{equation}
Note how the initial conditions $w(0,a)=1,\ w_x(0,a)=0$ imply $\rho=1,\ \vartheta=0$ at $x=0$.  
The first equation of \eqref{eq514} then implies $\rho>0$ on the rectangle $\mathcal{R}$.
In particular all zeros of $x\mapsto w(x,a)$ are simple, for any fixed $a$.
The second equation implies $\vartheta_x>0$ for $\cos \vartheta=0$, i.e. whenever $w=0$, alias 
$\vartheta=\frac12{\pi}\,(\mathrm{mod}\,\pi)$. 
Hence the total number of zeros of the solutions $x\mapsto w=w(x,a)$ relates to the winding of 
$x\mapsto \vartheta=\vartheta(x,a)$.

\begin{figure}[t] 
\centering \includegraphics[scale = 1.0]{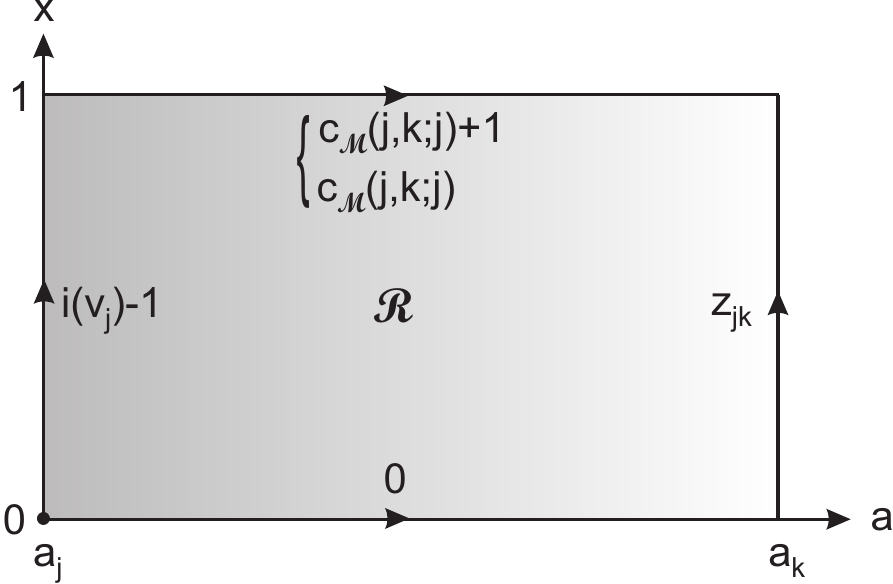}
\caption{\em\small Illustration of the rectangle $\mathcal{R}:=[0,1]\times[a_j,a_k]$. 
On each side of $\mathcal{R}$ is indicated the total winding of $\vartheta$ along the 
corresponding side. On the top side, the two values correspond to the alternative: $\cM_{j\,j+1}$ 
is in an odd/even quadrant with respect to $j$.} 
\label{fig51} 
\end{figure}

We can now outline the remaining proof of the NSL property \eqref{eq501} as follows.
We first note that the winding number of
\begin{equation}
\label{eq5wdg}
\vartheta\,(\mathrm{mod}\,2\pi): \quad \partial \mathcal{R} \rightarrow S^1
\end{equation}
is zero, along the boundary $\partial\mathcal{R}$ of the rectangle $\mathcal{R}$ in \eqref{eq5R}.
Indeed, the map extends to all of $\mathcal{R}$, continuously, and hence is contractible, i.e. of 
winding number zero.
In Lemma~\ref{l51} we relate the winding along the right boundary $a=a_k,\ 0\leq x\leq 1$ of 
$\mathcal{R}$ to the zero number $z(v_k-v_j)$ in claim \eqref{eq501}.
In Lemma~\ref{l52} we relate the winding along the left boundary $a=a_j,\ 0\leq x\leq 1$ of 
$\mathcal{R}$ to the Morse index $i(v_j)$ in claim \eqref{eq501}.
In Lemma~\ref{l53} we relate the winding along the upper boundary $x=1,\ a_j\leq a\leq a_k$ of 
$\mathcal{R}$ to the crossing number $c(j,k;j)$ in claim \eqref{eq501}.
Since $\vartheta=0$ is constant along the lower boundary $x=0,\ a_j\leq a\leq a_k$\,, and since the 
total winding number is zero, this reduces the proof of the NSL property \eqref{eq501} to the three 
Lemmata \ref{l51} -- \ref{l53}.

\begin{lemma} \label{l51}
The zero number $z(v_k-v_j)$ is given by
\begin{equation}
\label{eq515}
z(v_k-v_j)=z(w(\cdot,a_k))=\vartheta(1,a_k)/\pi\,.
\end{equation}
\end{lemma}

\emph{Proof:}
Fix $a=a_k$. We recall that all zeros of $w(\cdot,a_k)$ are simple.
They correspond to clockwise crossings of $(w,w_x)$ through the vertical $w_x$-axis or, equivalently, to 
simple zeros of $\vartheta-\tfrac{1}{2}\pi\,(\mathrm{mod}\,\pi)$ with positive slope $\vartheta_x>0$\,. 
The Neumann boundary condition $v_{k\,x}=v_{j\,x}=0$ at $x=1$ implies $w_x=0$ and hence 
$\vartheta\equiv 0\,(\mathrm{mod}\,\pi)$ there.
This proves the lemma. \hfill $\square$

\begin{lemma} \label{l52}
The Morse index $i(v_j)$, i.e.~the unstable dimension of the equilibrium $v_j$, satisfies
\begin{equation} \label{eq516}
i(v_j) = \lfloor \vartheta(1,a_j)/\pi \rfloor+1\,. 
\end{equation}
Here $\lfloor\cdot\rfloor$ denotes the integer valued floor function. 
\end{lemma}

\emph{Proof:} 
See for example \cite[Theorem 2]{roc85}.

Comparing \eqref{eq502} with \eqref{eq5w}, \eqref{eq507}, \eqref{eq509}, we first note that $w=v_a$ is the 
partial derivative of $v=v(\cdot,a)$ with respect to $a$, at $a=a_j$.
In particular, $(w(1,a_j),w_x(1,a_j))$ is the tangent of the meander segment $\cM_{j\,j+1}$ at $j$, at 
(clockwise) angle $\vartheta(1,a_j)\,(\mathrm{mod}\,2\pi)$ from the horizontal axis; see \eqref{eq5theta}.
Since the equilibrium $v_j$ is assumed to be hyperbolic, $\vartheta(1,a_j) \neq 0 \,(\mathrm{mod}\,\pi)$.
In particular the meander crosses the horizontal axis transversely at the intersection $j$.
More precisely, the (clockwise) tangent angle $\vartheta(1,a_j)$ has to point above the horizontal axis, for 
odd $j$ (alias even $i(v_j)$), and below for even $j$ (alias odd $i(v_j)$), alternatingly:
\begin{equation}
\label{eq517c}
\vartheta(1,a_j) \,(\mathrm{mod}\,2\pi)\in\begin{cases}
(0,\pi)\,,& \text{ for odd }i(v_j)\,; \\  
(\pi,2\pi)\,,& \text{ for even }i(v_j)\,.
\end{cases}
\end{equation}

Next, consider the simple eigenvalues $\lambda_m$\,, i.e.
\begin{equation} \label{eq511}
\lambda_0 > \dots > \lambda_{i(v_j)-1} > 0 > \lambda_{i(v_j)} > \dots \,,
\end{equation}
of the linearization 
\begin{equation} \label{eq510}
\lambda \tilde v=\tilde v_{xx}+q_1(x,a_j)\tilde v_x+q_0(x,a_j)\tilde v \,, 
\quad \tilde v(0)=1, \quad \tilde v_x(0)=\tilde v_x(1)=0
\end{equation}
at $v_j$. 
By classical Sturm-Liouville theory, e.g. as in \cite{cole}, the eigenfunction $\tilde v = \tilde v_m$ of 
$\lambda=\lambda_m$ possesses $m$ simple zeros. 
Moreover, $w$ solves \eqref{eq510} with $\lambda=0$, but violates the Neumann boundary condition at $x=1$. 
Therefore Sturm-Liouville comparison with 
\eqref{eq511} implies 
\begin{equation}
\label{eq517a}
i(v_j)-1 \leq z(w) \leq i(v_j)\,.
\end{equation}
Translating zeros of $w$ to zeros of $\vartheta - \tfrac{1}{2} \pi\,(\mathrm{mod}\,\pi)$, as in the proof of 
Lemma~\ref{l51}, we obtain
\begin{equation}
\label{eq517b}
(i(v_j)-1)\pi \ <\ \vartheta(1,a_j)- \pi/2 \ <\ (i(v_j)+1)\pi\,.
\end{equation}
Combined with \eqref{eq517c}, this implies
\begin{equation}
\label{eq517d}
(i(v_j)-1)\pi \ <\ \vartheta(1,a_j) \ <\ i(v_j)\pi
\end{equation}
and proves claim \eqref{eq516} of the lemma. \hfill $\square$

\medskip

We recall the definition of the signed clockwise counts $c(j,k;j)$ of meander crossings, 
from section \ref{sec2}.

\begin{lemma} \label{l53}
The clockwise increase 
\begin{equation}
\label{eq519}
\lfloor(\vartheta(1,a_k)-\vartheta(1,a_j))/\pi\rfloor=\vartheta(1,a_k)/\pi-\lfloor\vartheta(1,a_j))/\pi\rfloor
\end{equation} 
of the angle $\vartheta(1,a)$ from $a=a_j$ to $a=a_k$ is given by
\begin{equation} \label{eq512}
\begin{cases}
\  c(j,k;j)+1 , & \text{ if }\ Q\cM_{j\,j+1} \text{ is odd } \,; \\  
\\
\  c(j,k;j), & \text{ if }\ Q\cM_{j\,j+1} \text{ is even } \,.
\end{cases}
\end{equation}
\end{lemma}

\emph{Proof:}
At $a=a_k$ we recall $\vartheta(1,a)/\pi\in\N_0$ from \eqref{eq515}.  This proves claim \eqref{eq519}. 
It remains to prove claim \eqref{eq512}.

From the proof of Lemma~\ref{l52} we recall that $\vartheta(1,a_j)\,(\mathrm{mod}\,2\pi)$ is the (clockwise) 
tangent angle of the meander segment $\cM_{j\,j+1}$ at $j$ with the horizontal axis; see \eqref{eq5theta}.
For general $a_j<a\leq a_k$, the angle $\vartheta(1,a)$ tracks the secant between $j$ and 
$(w(1,a),w_x(1,a))\in \cM_{j\,k}$\,.

By definition, neither the clockwise crossing count $c(j,k;j)$, nor the clockwise angular increase 
$\lfloor(\vartheta(1,a_k)-\vartheta(1,a_j))/\pi\rfloor$ depend on homotopies of the Jordan meander segment 
$\cM_{jk}$, as long as the initial tangent angle $\vartheta(1,a_j)$ and the final secant $\vartheta(1,a_k)$ 
remain fixed.
Therefore we may assume all crossings of the angle $\vartheta(1,a)\,(\mathrm{mod}\,\pi)$ through the levels 
$\pi/2$ to be transverse, and finite in number, for $a\in(a_j,a_k]$.

The two cases \eqref{eq512} at $a=a_j$ arise as follows. 
The even/odd parity of $j$ at the up- or down-crossing $j$ of the meander $\cM$ determines whether the 
initial (clockwise) tangent $\vartheta(1,a_j)$ points down or up; see \eqref{eq517c}.
The precise direction, however, and the even/odd quadrant on that side of the horizontal axis, remain 
undetermined.
Suppose we artificially twist an initial tangent $\vartheta(1,a_j)$ from an odd quadrant at $j$ to the even 
quadrant on the same side of the horizontal axis.
We do not require this twist to be realized by specific nonlinearities $f=f(x,v,v_x))$ in \eqref{eq502}; we 
just compensate our twist, locally, by a homotopy of the initial meander segment $\cM_{j\,j+1}$ near $j$.
Then our twist of the initial tangent contributes one additional clockwise crossing to the crossing count 
$c(j,k;j)$.
Since this modification leaves \eqref{eq512} invariant, we may assume the initial tangent $\vartheta(1,a_j)$ 
to be in an odd quadrant, without loss of generality, i.e.
\begin{equation}
\label{eq518}
\lfloor\vartheta(1,a_j))/\pi\rfloor+\tfrac{1}{2}<\vartheta(1,a_j)/\pi<\lfloor\vartheta(1,a_j))/\pi\rfloor+1\,.
\end{equation}

Therefore \eqref{eq519}, \eqref{eq518} ensure that $\lfloor(\vartheta(1,a_k)-\vartheta(1,a_j))/\pi\rfloor -1$ 
counts the net clockwise crossings of the angle $\vartheta(1,a)\,(\mathrm{mod}\,\pi)$ through the levels 
$\pi/2$, as $a$ increases from $a=a_j$ to $a=a_k$\,.
This coincides with the clockwise crossing count $c(j,k;j)$, by definition, and proves the lemma.
\hfill $\square$

Contractibility of the winding map \eqref{eq5wdg} allows us to identify the values of $\vartheta$ in lemmata 
\ref{l51} -- \ref{l53} as real numbers, not just $\mathrm{mod}\,2\pi$. 
The lemmata therefore combine to show the NSL property \eqref{eq501} as follows. We first express 
$z_{j\,k}=\vartheta(1,a_k)/\pi$ in \eqref{eq515} of Lemma~\ref{l51} via \eqref{eq519} of Lemma~\ref{l53}. 
We then substitute the two floor function expressions in \eqref{eq519} by \eqref{eq512} and \eqref{eq516}, 
respectively. This proves our original claim \eqref{eq501}.
We may therefore evaluate $\vartheta(1,a_k)$ in \eqref{eq515} via summation of \eqref{eq514} and 
\eqref{eq519}, to complete the proof of the NSL property \eqref{eq501}.

\section{Appendix: Meander suspensions}\label{sec6}

The {\it double cone suspension} $\widetilde\cA$ of a global attractor $\cA$ is the topological quotient of 
$\cA\times[0,1]$ obtained by identifying $\cA\times 0$ and $\cA\times 1$ to distinct disjoint points, 
called {\it cone points}. 
The unstable suspension of global attractors is an efficient tool for the study of their geometric 
properties (see for example \cite{firo00}), and has also been considered in the related setting of 
meanders (see \cite{kar17}).
In this Appendix we define the corresponding suspension of a meander $\cM$ and review some results 
encoded in the boundary value orderings of the equilibria in their Sturm global attractors. 

\begin{figure} 
\centering \includegraphics[scale = 0.8]{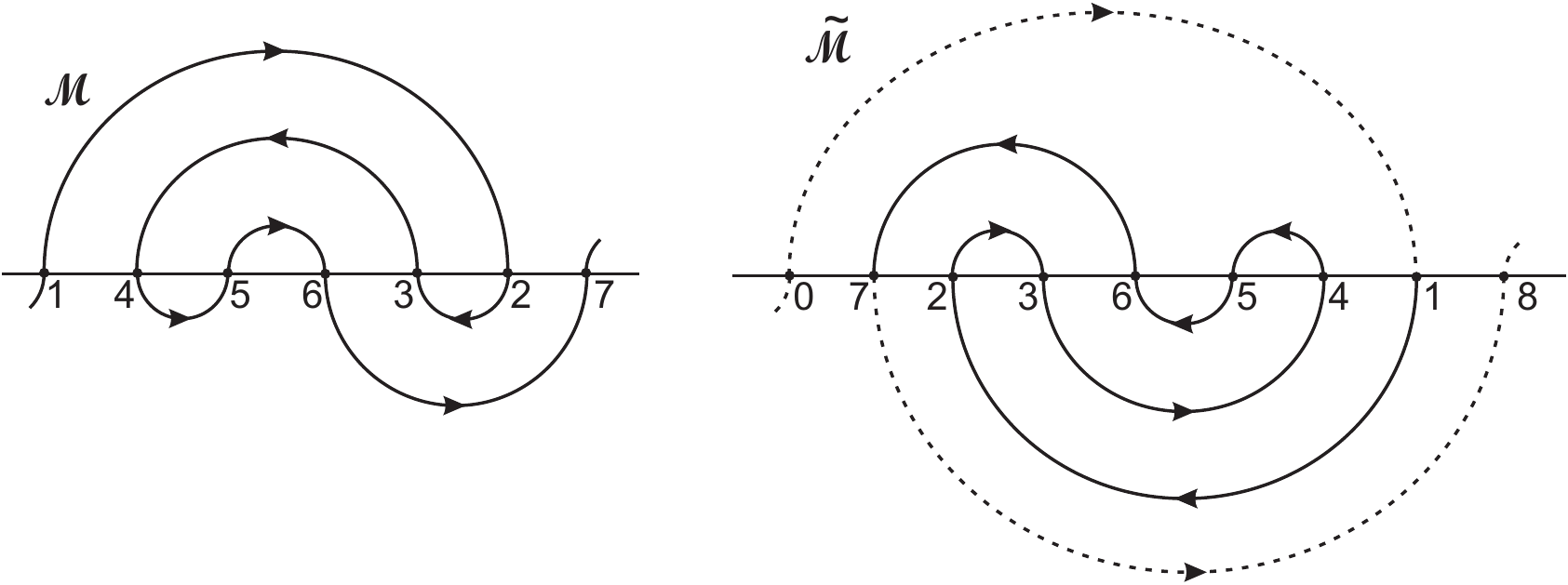}
\caption{\em\small Canonical form of the suspension $\widetilde\cM$ of the meander $\cM$ corresponding 
to the permutation $\sigma = \{ 1 \ 4 \ 5 \ 6 \ 3 \ 2 \ 7 \}$. The meander segment $\widetilde\cM_{1\,7}$ 
corresponds to a $180^\circ$ rotation of the segment $\cM_{1\,7}$. The suspended meander $\widetilde\cM$ 
corresponds to the permutation $\widetilde\sigma = \{ 0 \ 7 \ 2 \ 3 \ 6 \ 5 \ 4 \ 1 \ 8 \}$.}
\label{fig61}
\end{figure}

Let $\widetilde\cM$ denote the {\it suspension} of the meander $\cM$ obtained by rotating the segment 
$\cM_{1\,N}$ by $180^\circ$ and adding two extreme intersection points, labeled by $j=0$, left, 
and $j=N+1$, right, and two maximal outermost arcs. See Figure~\ref{fig61} for an illustration.

As we show in Lemma~\ref{l61} below, the suspension $\widetilde\cM$ of a Sturm meander is again Sturm 
and the corresponding global attractor $\widetilde\cA$ is connection equivalent to a double cone unstable 
suspension of the global attractor $\cA$. 
Let $\widetilde\cE$ denote the set of equilibria of $\widetilde\cA$ and 
$\widetilde h_0,\widetilde h_1:\{0,\dots,N+1\}\rightarrow\widetilde\cE$ their boundary orders as 
obtained from the suspension $\widetilde\cM$. Due to the rotation of the meander segment $\cM_{1\,N}$, 
the correspondence $h_0(j)=v_j\mapsto\widetilde v_j:=\widetilde h_0(j), j=1,\dots,N$, 
preserves the $h_0$-order and reverses the $h_1$. Moreover, as we will show, the first equilibrium 
$\widetilde v_0:=\widetilde h_0(0)$ and the last equilibrium $\widetilde v_{N+1}:=\widetilde h_0(N+1)$ 
constitute the cone points of $\widetilde\cA$. 

\begin{lemma} \label{l61}
Let $\widetilde\cM$ denote the meander suspension of the Sturm meander $\cM$. 
Then, $\widetilde\cM$ is a Sturm meander. Moreover, the Morse indices of the equilibria 
in $\widetilde\cE$ and their zero number relations satisfy  
\begin{equation} \label{eq601}
i(\widetilde v_0)=i(\widetilde v_{N+1})=0 \ , \quad 
i(\widetilde v_j) = i(v_{j})+1 \mbox{ for } \ 1\le j\le N \ ,
\end{equation}
\begin{equation} \label{eq602}
z(\widetilde v_j-\widetilde v_k)=z(v_j-v_k)+1 \ \mbox{ for } \ 1\le j<k\le N , 
\end{equation}
\begin{equation} \label{eq603}
z(\widetilde v_0-\widetilde v_j)=0 \ \mbox{ for } \ 1\le j\le N+1 , \quad
z(\widetilde v_k-\widetilde v_{N+1})=0 \ \mbox{ for } \ 0\le k\le N \,.
\end{equation}
\end{lemma}

\emph{Proof:}
Let $\sigma:=h_0^{-1}\circ h_1$ and $\widetilde\sigma:=\widetilde h_0^{-1}\circ\widetilde h_1$ denote 
the permutations corresponding to the Sturm meanders $\cM$ and $\widetilde\cM$, respectively.
Let $\widetilde\sigma'$ denote the restriction of $\widetilde\sigma$ to the set $\{1,\dots,N\}$. 
By our definition of meander suspension, the meander segment $\widetilde\cM_{1\,N}$ corresponds to 
a $180^\circ$ rotation of $\cM_{1\,N}$. This implies that $\widetilde\sigma'=\kappa\sigma$, where $\kappa$ 
is the reversal involution. Therefore, we have
\begin{equation} \label{eq604}
\widetilde\sigma(0)=0 \ , \quad \widetilde\sigma(N+1)=N+1 \ , \quad 
\widetilde\sigma(j)=\sigma(N+1-j) \ \mbox{ for } 1\le j\le N \,, 
\end{equation}
\begin{equation} \label{eq605}
\widetilde\sigma^{-1}(j)=\sigma^{-1}\kappa(j)=N+1-\sigma^{-1}(j) \ \mbox{ for } 1\le j\le N \,. 
\end{equation}
To show that $\widetilde\sigma\in S_{N+2}$ is a Sturm permutation, we compute the corresponding Morse 
numbers $\widetilde i_j, j=0,\dots,N+1$, using the recursion \eqref{eq201} in this setting:
\begin{equation} \label{eq606}
\widetilde i_0=0 \ , \quad \widetilde i_{j+1} = \widetilde i_j + 
(-1)^j\sign\left(\widetilde\sigma^{-1}(j+1)-\widetilde\sigma^{-1}(j)\right) \ \mbox{ for } 0\le j\le N \,. 
\end{equation}
In terms of the permutation $\sigma$, by \eqref{eq605}, this recursion becomes 
\begin{equation} \label{eq607}
\widetilde i_{j+1} = \widetilde i_j +  
(-1)^j\sign\left(\sigma^{-1}(j)-\sigma^{-1}(j+1)\right) \ \mbox{ for } 1\le j<N \,.
\end{equation}
A comparison between \eqref{eq607} and \eqref{eq201} then shows that both recursions, 
$\widetilde i_j$ and $i_j$, have the same step increments, 
\begin{equation} \label{eq608}
\widetilde i_{j+1}-\widetilde i_j = i_{j+1}-i_j \ \mbox{ for all } 1\le j<N \,.
\end{equation}
Now, these recursions start at $j=1$ with $i_1=0$ and $\widetilde i_1=1$, by \eqref{eq606}.   
Therefore, \eqref{eq608} shows that $\widetilde i_j=i_j+1$ for all $1\le j\le N$. 
It follows that $\widetilde i_N=1$, and \eqref{eq606} again implies $\widetilde i_{N+1}=0$. 
This shows that $\widetilde\cM$ is a Sturm meander and proves \eqref{eq601}.

Next, to prove the zero number relations \eqref{eq602}--\eqref{eq603} we invoke \eqref{eq204}, 
alias the nonlinear Sturm-Liouville property \eqref{eq501} in the appendix section~\ref{sec5}.

Let $c=c(j,k;\ell)$ and $\tilde c=\tilde c(j,k;\ell)$ denote the crossing numbers of $\cM_{1\,N}$ and 
$\widetilde\cM_{0\,N+1}$, respectively. We first note that the suspension does not affect the meander 
orientation. In fact, both meander segments $\cM_{1\,N}$ and $\widetilde\cM_{1\,N}$ are oriented by the 
increasing labels along $h_0$ and $\widetilde h_0$, respectively.
Since $\cM_{1\,N}$ and $\widetilde\cM_{1\,N}$ are equal up to rotation, all the crossing numbers are 
preserved by the meander suspension, i.e. we have
\begin{equation} \label{eq609}
\widetilde c(j,k;\ell) = c(j,k;\ell) \ \mbox{ for all } 1 \le j,k,\ell\le N \,. 
\end{equation}
The zero number relations for the equilibria $\widetilde v_1,\dots,\widetilde v_N\in\widetilde\cE$ 
are then obtained directly from \eqref{eq204}. For all $0\le j<k\le N+1$, we have
\begin{equation} \label{eq610}
z(\widetilde v_k-\widetilde v_j) = 
\begin{cases}
i(\widetilde v_j) + \widetilde c(j,k;j), & \text{ if } 
Q\widetilde \cM_{j\,j+1} \text{ is odd } ; \\  \\
i(\widetilde v_j) - 1 + \widetilde c(j,k;j), & \text{ if } 
Q\widetilde \cM_{j\,j+1} \text{ is even } . 
\end{cases}
\end{equation}
In addition, for $1\le j<N$ the quadrants $Q\widetilde \cM_{j\,j+1}$ and $Q\cM_{j\,j+1}$ have the same 
even/odd parity. Indeed, the parity of the quadrants is not affected by the $180^\circ$ rotation. 
Then, by \eqref{eq601} and \eqref{eq609}, from \eqref{eq610} we obtain for all $1\le j<k\le N$ 
\begin{equation} \label{eq611}
z(\widetilde v_k-\widetilde v_j) = 
\begin{cases}
i(v_j) + 1 + c(j,k;j), & \text{ if } Q\cM_{j\,j+1} \text{ is odd } ; \\  \\
i(v_j) + c(j,k;j), & \text{ if } Q\cM_{j\,j+1} \text{ is even } . 
\end{cases}
\end{equation}
Hence, comparing \eqref{eq611} with \eqref{eq204} for the equilibria $v_1,\dots,v_N\in\cE$, 
we obtain \eqref{eq602}.

Finally, to prove \eqref{eq603} we invoke again \eqref{eq204} for the suspended meander $\widetilde\cM$, 
alias \eqref{eq610}. By the definition of meander suspension, the quadrant $Q\widetilde\cM_{0\,1}$ 
of the first arc segment of $\widetilde\cM$ is odd. 
Moreover, $i(\widetilde v_0)=0$ for the extremal equilibrium $\widetilde v_0$, and the crossing numbers 
of $\widetilde\cM_{0\,k}$ with respect to $\ell=0$ satisfy 
\begin{equation} \label{eq612}
\widetilde c(0,k;0) = 0 \ \mbox{ for all } 1 \le k\le N+1 \,. 
\end{equation}
Therefore, \eqref{eq610} with $j=0$ implies
\begin{equation} \label{eq613}
z(\widetilde v_k-\widetilde v_0) = i(\widetilde v_0) + \widetilde c(0,k;0) = 0 \ 
\mbox{ for all } 1\le k\le N+1 \,.
\end{equation}
Similarly, for the extremal equilibrium $\widetilde v_{N+1}$, we obtain
\begin{equation} \label{eq614}
z(\widetilde v_{N+1}-\widetilde v_j) = 0 \ \mbox{ for all } 0\le j\le N \,.
\end{equation}
This proves \eqref{eq603}, and completes the lemma.
\hfill $\square$

\medskip

Let $\cE^*:=\widetilde\cE\setminus\{\widetilde v_0,\widetilde v_{N+1}\}$ and let $\cH^*$ denote 
the set of heteroclinic orbits between the equilibria in $\cE^*$, i.e. the connecting orbits $v\leadsto w$ 
with $v,w\in\cE^*$. Let $\cA^*\subset\widetilde\cA$ denote the subset of all equilibria in $\cE^*$ and 
all their connecting orbits,
\begin{equation} \label{eq615}
\cA^* := \cE^* \cup \cH^* \,.
\end{equation}
Then, $\cA^*$ is a flow invariant subset of $\widetilde\cA$ which, by the transversality of stable and 
unstable manifolds, is the union of invariant submanifolds of $\widetilde\cA$. 

By our previous Lemma~\ref{l61}, the suspension has the effect of increasing all Morse indices by one, 
as asserted by \eqref{eq601}. Moreover, the uniform increase by one of all the 
zero number relations between the equilibria, displayed by \eqref{eq602}, shows that the suspension 
preserves the heteroclinic structure of the global attractor. Hence, the suspension correspondence 
$\cE\mapsto\cE^*$ determined by $v_j\mapsto\widetilde v_j, j=1,\dots,N$, 
implies the connection equivalence between the global attractor $\cA$ and the subset $\cA^*$, i.e. 
the corresponding oriented connection graphs are isomorphic,
\begin{equation} \label{eq616}
\cA \sim \cA^* \,.
\end{equation}

\begin{figure}
\centering \includegraphics[scale = 1.2]{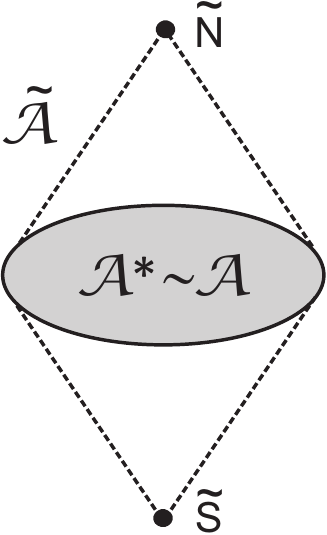}
\caption{\em\small Illustration of the unstable suspension of the global attractor $\cA$. 
The cones are topologically glued at their bases where they share an invariant subset $\cA^*$ which is 
connection equivalent to the global attractor $\cA$, $\cA^*\sim\cA$. The cone points $\bf\widetilde N$ 
and $\bf\widetilde S$ correspond to the extremal equilibria $\widetilde v_0$ and $\widetilde v_{N+1}$.}
\label{fig62}
\end{figure}

The zero number relations \eqref{eq603} then concern the orbits contained in $\widetilde\cA\setminus\cA^*$. 
In particular, for each $\widetilde v_j\in \cA^*$, \eqref{eq603} and the parabolic comparison principle 
imply that the one-dimensional fast unstable manifold $W^1(\widetilde v_j)$ 
is the union of $\widetilde v_j$ itself, with two heteroclinic orbits: one  
connecting to $\widetilde v_0$, $\widetilde v_j\leadsto\widetilde v_0$, and the other connecting to 
$\widetilde v_{N+1}$, $\widetilde v_j\leadsto\widetilde v_{N+1}$. Therefore, $\cA^*$ stands unstably 
suspended between the two extremal equilibria, $\widetilde v_0$ and $\widetilde v_{N+1}$, introduced by 
the suspension. See Figure~\ref{fig62} for an illustration.

The next Lemma~\ref{l62} addresses the behavior of the minimax property under the meander suspension.

\begin{lemma} \label{l62}
Let $\widetilde\cO\in\widetilde\cE\setminus\{\widetilde v_0,\widetilde v_{N+1}\}$ denote the equilibrium 
corresponding to any unstable $\cO\in\cE$ by the suspension mapping $v_j\mapsto\widetilde v_j, j=1,\dots,N$.
Then we have the following correspondences between neighbors and minimax equilibria:
\begin{equation} \label{eq617}
\cE^{n-1}_\pm(\cO) \mapsto \cE^n_\pm(\widetilde\cO) \ , \quad 
\underline v^{\,\iota}_{\,\pm}(\cO) \mapsto \underline v^{\,\iota}_{\,\pm}(\widetilde\cO) \ , \ 
\overline v^{\,\iota}_{\,\pm}(\cO) \mapsto \overline v^{\,\iota}_{\,\pm}(\widetilde\cO) \ , \ 
\iota\in\{0,1\} \,.
\end{equation}
\end{lemma}

\emph{Proof:}
The suspended equilibria correspondence $\cE^{n-1}_\pm(\cO) \mapsto \cE^n_\pm(\widetilde\cO)$, 
follows directly from \eqref{eq603} of Lemma~\ref{l61} and the connection equivalence $\cA\sim\cA^*$, 
since the meander suspension preserves the $h_0$-order.
The same arguments prove the correspondences 
$\underline v^{\,\iota}_{\,\pm}(\cO)\mapsto\underline v^{\,\iota}_{\,\pm}(\widetilde\cO)$ and 
$\overline v^{\,\iota}_{\,\pm}(\cO)\mapsto\overline v^{\,\iota}_{\,\pm}(\widetilde\cO)$ for 
$\iota\in\{0,1\}$, and complete the proof of \eqref{eq617} and the lemma.
\hfill $\square$

\medskip

In conclusion, this lemma shows that the Morse index $n=i(\cO)$ can indeed be assumed odd, without any loss, 
in our proof of Theorem~\ref{theo11}.
In fact, if the minimax property $\underline v^\iota_\pm(\cO)=\overline v^{1-\iota}_\pm(\cO)$ holds for 
odd Morse index $i(\cO)$, Lemma~\ref{l62} asserts that the minimax property 
$\underline v^{\,\iota}_{\,\pm}(\widetilde\cO)=\overline v^{\,1-\iota}_{\,\pm}(\widetilde\cO)$ holds 
for even Morse index $i(\widetilde\cO)=i(\cO)+1$. 

As a final observation we mention that the $180^\circ$ rotation introduced in our definition of meander 
suspension is not necessary. In fact, Lemma~\ref{l61} and also Lemma~\ref{l62} conveniently adapted, 
hold as well for a meander suspension defined without the rotation of the meander segment $\cM_{1\,N}$.


\begin{thebibliography}{20}

\bibitem[An88]{ang88}
S.~Angenent.
The zero set of a solution of a parabolic equation.
{\em J.~Reine Angew. Math.}, \textbf{390} (1988), 79--96.

\bibitem[Ar88]{arn88}
V.~I.~Arnold. 
A branched covering $CP^2\rightarrow S^4$, hyperbolicity and projective topology.
{\em Siberian Math. J.}, \textbf{29} (1988), 717--726.

\bibitem[BV92]{bavi92}
A.~V.~Babin and M.~I.~Vishik.
{\em Attractors of Evolution Equations}.
North Holland, Amsterdam, 1992.   

\bibitem[BF88]{brfi88}
P.~Brunovsk\'y and B.~Fiedler.
Connecting orbits in scalar reaction diffusion equations.
{\em Dynamics Reported} \textbf{1} (1988), 57--89.

\bibitem[BF89]{brfi89}
P.~Brunovsk\'y and B.~Fiedler.
Connecting orbits in scalar reaction diffusion equations {II}: The complete solution.
{\em J.~Differential~Eqs.}, \textbf{81} (1989), 107--135.

\bibitem[CL72]{cole}
E.A.~Coddington and N.~Levinson.
\emph{Theory of Ordinary Differential Equations.}
McGraw-Hill, New York 1974.

\bibitem[D\&al19]{detal19} 
V.~Delecroix, \'E.~Goujard, P.~Zograf, and A.~Zorich. 
Enumeration of meanders and Masur-Veech volumes. 
arxiv:1705.05190v2; Forum of Mathematics, Pi \textbf{8} (2020), 80pp;
doi: 10.1017/fmp.2020.2

\bibitem[Fi94]{fie94}
B.~Fiedler.
Global attractors of one-dimensional parabolic equations: Sixteen examples.
{\em Tatra Mountains Math.~Publ.} \textbf{4} (1994), 67--92.

\bibitem[FR96]{firo96} 
B.~Fiedler and C.~Rocha.  
Heteroclinic orbits of semilinear parabolic equations. 
{\em  J.~Differential~Eqs.}, \textbf{125} (1996), 239--281. 

\bibitem[FR99]{firo99} 
B.~Fiedler and C.~Rocha.
Realization of meander permutations by boundary value problems.  
{\em J.~Differential~Eqs.}, \textbf{156} (1999), 282--308.

\bibitem[FR00]{firo00}
B.~Fiedler and C.~Rocha.
Orbit equivalence of global attractors of semilinear parabolic differential equations.
{\em Trans.~Amer.~Math.~Soc.}, \textbf{352} (2000), 257--284.

\bibitem[FR08]{firo08} 
B.~Fiedler and C.~Rocha.
Connectivity and design of planar global attractors of Sturm type.
II: Connection graphs. 
{\em J.~Differential~Eqs.}, \textbf{244} (2008), 1255--1286. 

\bibitem[FR09a]{firo09a} 
B.~Fiedler and C.~Rocha.
Connectivity and design of planar global attractors of Sturm type. 
I: Orientations and Hamiltonian paths. 
{\em J.~Reine Angew. Math.}, \textbf{635} (2009), 71--96. 

\bibitem[FR09b]{firo09b} 
B.~Fiedler and C.~Rocha.
Connectivity and design of planar global attractors of Sturm type.
III: Small and Platonic examples. 
{\em J.~Dynam. Differential Eqs.}, \textbf{22} (2010), 121--162.

\bibitem[FR14]{firo14} 
B.~Fiedler and C.~Rocha. 
Nonlinear Sturm global attractors: unstable manifold decompositions as regular CW-complexes. 
{\em Discrete Contin.~Dyn.~Syst.}, \textbf{34} (2014), no. 12, 5099–-5122.

\bibitem[FR15]{firo15} 
B.~Fiedler and C.~Rocha. 
Schoenflies spheres as boundaries of bounded unstable manifolds in gradient Sturm systems. 
{\em Discrete Contin.~Dyn.~Syst.}, \textbf{27} (2015), no. 3-4, 597–-626.

\bibitem[FR18a]{firo18a} 
B.~Fiedler and C.~Rocha. 
Sturm 3-ball global attractors 1: Thom-Smale complexes and meanders. 
{\em S\~ao~Paulo~J.~Math.~Sci.}, \textbf{12} (2018), 18--67.

\bibitem[FR18b]{firo18b} 
B.~Fiedler and C.~Rocha. 
Sturm 3-ball global attractors 2: Design of Thom-Smale complexes. 
{\em J.~Dynam.~Differential Eqs.}, \textbf{31} (2019), 1549--1590.

\bibitem[FR18c]{firo18c} 
B.~Fiedler and C.~Rocha. 
Sturm 3-ball global attractors 3: Examples of Thom-Smale complexes. 
{\em Discrete Contin.~Dyn.~Syst.}, \textbf{38} (2018), no. 7, 3479--3545.

\bibitem[FR20]{firo18d} 
B.~Fiedler and C.~Rocha.
Boundary orders and geometry of the signed Thom-Smale complex for Sturm global attractors.
arxiv: 1811.04206; \emph{J.~Dynam.~Differential Eqs.} (2020), 32pp.; doi: 10.1007/s10884-020-09836-5

\bibitem[FRW12]{firowo12} 
B.~Fiedler, C.~Rocha, M.~Wolfrum:
A permutation characterization of Sturm global attractors of Hamilton type.
\emph{J.~Differential Eqs.} \textbf{252} (2012), 588--623.

\bibitem[FR91]{furo91}
G.~Fusco and C.~Rocha.
A permutation related to the dynamics of a scalar parabolic {PDE}.
{\em J.~Differential Eqs.}, \textbf{91} (1991), 75--94.

\bibitem[Ga04]{gal04}
V.~A.~Galaktionov.
{\em Geometric Sturmian Theory of Nonlinear Parabolic Equations and Applications}.
Chapman \& Hall, Boca Raton, 2004.

\bibitem[Ha88]{hal88}
J.~K.~Hale.
{\em Asymptotic Behavior of Dissipative Systems}. 
Math.~Surv. \textbf{25}. AMS Publications, Providence, 1988.

\bibitem[He81]{hen81}
D.~Henry.
{\em Geometric Theory of Semilinear Parabolic Equations}.
Lect.~Notes~Math. \textbf{804}, Springer-Verlag, New York, 1981.

\bibitem[Ka17]{kar17}
A.~Karnauhova.
{\em Meanders -- Sturm Global Attractors, Seaweed Lie Algebras and Classical Yang-Baxter Equation}.
De Gruyter, Berlin, 2017.

\bibitem[La91]{lad91}
O.~A.~Ladyzhenskaya.
{\em Attractors for Semigroups and Evolution Equations}.
Cambridge University Press, 1991.

\bibitem[Ma82]{mat82}
H.~Matano.
Nonincrease of the lap-number of a solution for a one-dimensional semi-linear parabolic equation.
{\em J.~Fac.~Sci.~Univ.~Tokyo Sec.~IA.}, \textbf{29} (1982), 401--441.

\bibitem[Pa83]{paz83}
A.~Pazy.
{\em Semigroups of Linear Operators and Applications to Partial Differential Equations}.
Springer-Verlag, New York, 1983.

\bibitem[Ra02]{rau02}
G.~Raugel.
Global attractors in partial differential equations.
{\em Handbook of Dynamical Systems}, \textbf{2} (2002), 885--982.

\bibitem[Ro85]{roc85} 
C.~Rocha. 
Generic properties of equilibria of reaction-diffusion equations with variable diffusion. 
{\em Proc. Roy. Soc. Edinburgh Sect. A}, \textbf{101} (1985), 385--405.

\bibitem[Ro91]{roc91} 
C.~Rocha. 
Properties of the attractor of a scalar parabolic PDE. 
{\em J.~Dynam.~Differential Eqs.}, \textbf{3} (1991), 575--591.

\bibitem[Ta79]{tan79}
H.~Tanabe.
{\em Equations of Evolution}.
Pitman, Boston, 1979.

\bibitem[Te88]{tem88}
R.~Temam.
{\em Infinite-Dimensional Dynamical Systems in Mechanics and Physics}.
Springer-Verlag, New York, 1988.

\bibitem[Wo02]{wol02}
M.~Wolfrum.
A sequence of order relations: encoding heteroclinic connections in scalar parabolic {PDE}.
{\em J. Differential Eqs.} \textbf{183} (2002), no. 1, 56--78.

\end{thebibliography}
\end{document}